 \documentclass{siamart1116}
 
\usepackage{lipsum}
\usepackage{amsfonts}
\usepackage{graphicx}
\usepackage{epstopdf}
\usepackage{algorithmic}
\ifpdf
  \DeclareGraphicsExtensions{.eps,.pdf,.png,.jpg}
\else
  \DeclareGraphicsExtensions{.eps}
\fi

\newcommand{\TheTitle}{A Krylov-Schur like method for computing the
  best  rank-{$(r_1,r_2,r_3)$}  approximation of large and sparse tensors}

\title{{\TheTitle}}

\author{%
  Lars Eld\'en\thanks{Department of Mathematics, Link\"{o}ping
    University, SE-581 83, Link\"{o}ping, Sweden
    (\email{lars.elden@liu.se})}%
  \and
 Maryam Dehghan\thanks{Department of Mathematics, Persian Gulf University, 75169, Bushehr, Iran (\email{ma.dehghan@mehr.pgu.ac.ir, maryamdehghan880@yahoo.com}),}
}

\usepackage{amsopn}

\DeclareMathOperator{\gap}{gap}  %
\newcommand{\cA}{\mathcal{A}}
\newcommand{\cB}{\mathcal{B}}
\newcommand{\cC}{\mathcal{C}}
\newcommand{\cD}{\mathcal{D}}

\newcommand{\cF}{\mathcal{F}}

\newcommand{\cH}{\mathcal{H}}

\newcommand{\cU}{\mathcal{U}}
\newcommand{\cX}{\mathcal{X}}

\newcommand{\cW}{\mathcal{W}}

\newcommand{\RR}{\mathbb{R}}

\newcommand{\Gr}{\mathrm{Gr}}

\renewcommand{\a}{\alpha}
\renewcommand{\b}{\beta}
\newcommand{\g}{\gamma}
\renewcommand{\l}{\lambda}

\newcommand{\<}{\langle}
\renewcommand{\>}{\rangle}

\newcommand{\x}{\times}

\newcommand{\norm}[1]{\Vert #1 \Vert}

\newcommand{\tp}{^{\sf T}}

\newcommand{\tml}[3][]{\bm{\left(} #2 \bm{\right)}_{ #1} \bm{\cdot} #3}
\newcommand{\tmr}[3][]{#2 \bm{\cdot} \bm{\left(} #3 \bm{\right)}_{ #1}}

\newcommand{\matlab}{\textsc{matlab}}

\newcommand{\svd}{\textsc{svd}}

  \usepackage{bm}
  \usepackage{hhline}
  \usepackage{url}
  \usepackage{algorithm,algorithmic}
 \ifpdf
 \fi 
  
\begin{document}

 \maketitle

 \begin{abstract}
 The paper is concerned with methods 
 for computing the best low multilinear rank approximation of large
 and sparse tensors. Krylov-type methods have been used for this
 problem; here block versions  are introduced. 
For the computation of partial eigenvalue and singular value
decompositions of matrices the Krylov-Schur (restarted Arnoldi) method
is used. We describe a generalization of this method to tensors, for
 computing the best low multilinear rank 
 approximation of large and sparse tensors. In analogy to the matrix
 case, the large  tensor is only accessed in multiplications between the
 tensor and blocks of vectors, thus avoiding
 excessive memory usage. It is  proved that, if the  starting
 approximation is  good enough, then the tensor Krylov-Schur method is
 convergent. Numerical examples are given for synthetic tensors  and
 sparse tensors from applications, which demonstrate that for most
 large problems  the Krylov-Schur method converges faster and more robustly
 than higher order orthogonal iteration. 
 \end{abstract}

 \begin{keywords}
   tensor, multilinear rank, best rank-(p,q,r) approximation,
   Grassmann manifold, sparse tensor, block Krylov-type method, Krylov-Schur
   algorithm, (1,2)-symmetric tensor
 \end{keywords}       

 \begin{AMS}
 65F99,  15A69, 65F15.
 \end{AMS}

 \section{Introduction}\label{sec:introd} In many applications of
 today,  large and sparse data 
 sets are generated that are organized in more than two
 categories. Such multi-mode data can be represented by tensors. They arise in
 applications of data sciences, such as web link 
 analysis, cross-language information
 retrieval and social network analysis, see, e.g.,
 \cite{kolda2005higher} and the survey \cite{kolda2009tensor}. The
 effective analysis of tensor data 
 requires the development of methods that can identify the inherent
 relations that exist in the data, and that scale to
 large data sets. Low rank approximation is one such method, and much
 research has been done in recent years in this area; a few examples
 that are related to the work of this paper are given in
 \cite{andersson1998improving, comon2002tensor, de2000best, 
   elden2009newton, ishteva2011best, khoromskij2009multigrid,
   oseledets2008tucker, savas2013krylov, savas2010quasi,
   zhang2001rank}. However, most of these methods are intended for
 small to medium size tensors. 
 The objective of this paper is to develop
  an algorithm for low-rank  approximation  of large and
 sparse tensors.

 We consider the problem of approximating a 3-mode tensor $\cA$  by another
 tensor $\cB$, 
 \begin{equation}\label{eq:min}
 \min\limits_{\cB}\norm{\cA-\cB},
 \end{equation}
 where the norm is the Frobenius norm, and $\cB$ has low multilinear
 rank-$(r_1,r_2,r_3)$ (for 
   definitions of the concepts used in this introduction, see 
    Section \ref{sec:preliminary}). We will assume that $\cA$ is large and
   sparse. This problem can be written
    \begin{equation}\label{eq:main min-intro}
 \min_{\cF,U,V,W}\norm{\cA-\tml{U,V,W}{\cF}},\quad
 \mathrm{subject~  to} \quad  U^TU=I_{r_1},\quad V^TV=I_{r_2},\quad W^TW=I_{r_3}, 
 \end{equation}
 where $\cF \in \RR^{r_1 \x r_2 \x r_3}$ is a tensor of  small
 dimensions, and $\tml{U,V,W}{\cF}$ 
 denotes matrix-tensor multiplication in all three modes.
This is the best rank-$(r_1,r_2,r_3)$  approximation problem
\cite{de2000best}, and it is a special case of  Tucker tensor
approximation \cite{tu:64,tu:66}. 
 It can be considered as a generalization of the problem of
computing the  Singular Value Decomposition (SVD) of a matrix
\cite{elden2011perturbation}. In fact, a partial SVD solves the
matrix approximation problem corresponding to \eqref{eq:min}, see
e.g. \cite[Chapter 2.4]{golub2013matrix}.   

In this paper we develop a Block Krylov-Schur like (BKS) method for
computing the best rank-$(r_1,r_2,r_3)$ approximation of large and
sparse tensors which are  symmetric with respect to two modes. We are
specially interested in small values of the rank, as  in the
two parallel papers \cite{eldehg20a,eldehg20c}, where  it is
essential to use the best approximation rather than any Tucker
approximation. 

 Krylov methods are routinely used to compute partial SVD's (and
 eigenvalue decompositions) of large sparse  matrices
 \cite{lehoucq1998arpack}. In \cite{savas2013krylov} we introduced a
 generalization of Krylov methods to tensors. 
 It was shown experimentally  that tensor
 Krylov-type methods have similar approximation  properties as the
 corresponding matrix Krylov methods.  In this paper we present
 \emph{block versions of Krylov-type methods for tensors}, which are expected
 to be more efficient than the methods in
 \cite{savas2013krylov}. Having problems in mind, where the tensor is
  symmetric with respect to  two modes (e.g. sequences of
 adjacency matrices of graphs), we formulate the block methods in
 terms of such tensors. 

 Even if matrix Krylov methods give low rank approximations, their
 convergence properties are usually not good enough, and, if used
 straightforwardly, they may require excessive memory and computer
 time. Therefore they are accelerated using  restart techniques
 \cite{lehoucq1998arpack}, which are equivalent to the  Krylov-Schur
 method \cite{stewart2002krylov}.  We here present a \emph{tensor
   Krylov-Schur like method},  and show that it can be used to compute best
 multilinear  low-rank approximations of large and sparse tensors.

 The Krylov-Schur method is an inner-outer iteration. In the outer
 iteration we start from an approximate solution and generate new
 blocks of orthogonal vectors, using Krylov-type block methods. Thus
 we use the large tensor only in tensor-matrix multiplications, where
 the matrix has few columns (this is analogous to the use of block Krylov
 methods for large matrices). Then
 we project the problem \eqref{eq:main min-intro} to a smaller problem
 of the same type, which we solve in the inner iterations, using a
 method for problems with a medium-size, dense tensor. The problem can
 be formulated as one on a product of Grassmann manifolds
 \cite{elden2009newton}. In our experiments we use a Newton-Grassmann
 method. As a stopping criterion for the outer  iteration we use the
 norm of the Grassmann gradient of the objective function. We show
 that this gradient can be computed efficiently in terms of the small
 projected problem. 

 We also prove that, if the starting approximation the Krylov-Schur
 method is good enough, the BKS method is convergent.

 The literature on algorithms for best rank-$(r_1,r_2,r_3)$
 approximation of large and sparse tensors is not extensive
 \cite{kosu08,goreinov2012wedderburn,kauc15}.  To the best of our
 knowledge the present paper is the first one that goes significantly beyond
the  Higher Order Orthogonal Iteration (HOOI) \cite{de2000best}. Our
experiments indicate that the new method is more efficient and robust than HOOI
for large and sparse tensors. 

  The paper is organized as follows. Some tensor
  concepts are introduced in Section \ref{sec:preliminary}. The Krylov-Schur
  procedure for matrices is  sketched in Section
  \ref{sec:krylov-type}.  In Section \ref{sec:Krylov-tensor} block
  Krylov-type methods for tensors are   described. The tensor
  Krylov-Schur method is presented and analyzed in Section
  \ref{sec:EK-Schur-gen}.       Some
  numerical examples are given in Section \ref{sec:numerical test}
  that illustrate the accuracy and efficiency 
  of the Krylov-Schur method. In particular, we demonstrate that for
  less well-conditioned approximation problems, the new method converges
  faster than HOOI \cite{de2000best}. 
 We are especially interested in tensors that are symmetric
  with respect to the first two modes (e.g. tensor consisting of a
  sequence of adjacency matrices of undirected graphs). Therefore the
  block Krylov-type methods and the numerical examples are given in
  terms of such tensors.

  The method presented in this paper is ``like'' a Krylov-Schur method for
  two reasons. The tensor Krylov-type method  is not a  Krylov
  method in a strict sense,  as it does not build bases for Krylov subspaces
  \cite{savas2013krylov}. The  method is not a real
  Krylov-Schur method as it does not build and manipulate a Hessenberg
  matrix; instead it uses a tensor, which is in some ways 
  similar to Hessenberg. However, this structure is not utilized.
  For ease of presentation we will sometimes omit
  ``like'' and ``type''.

  \section{Tensor concepts and preliminaries}\label{sec:preliminary}  
 \subsection{Notation}\label{subsec:notations}
 Throughout this paper we use of the following notations. Vectors will
 be denoted by lower case roman letters, e.g., $a$
 and $b$, matrices by capital roman letters, e.g., $A$ and $B$
 and tensors by calligraphic letters, e.g., ${\cA}$ and
 ${\cB}$.

  Notice that sometimes we will not explicitly mention the dimensions
  of matrices and tensors, and assume that they are such that the
  operations are well-defined. Also, for simplicity of notation and
  presentation, we will restricted ourselves to tensors of order 3,
 which are defined in the next paragraph. The
  generalization to higher order tensors is straightforward. For more
  general definitions we refer reader to \cite{bader2006algorithm}. 

  A tensor will be viewed as a multidimensional array of real
  numbers. The order of a tensor is the number of dimensions, also
  known as modes, e.g., a 3-dimensional array,
  $\cA\in \RR^{l\x m\x n}$, is called a tensor of order 3 or
  3-tensor. A fiber is a one-dimensional section of a tensor, obtained
  by fixing all indices except one; $\cA(i,:,k)$ is referred to
  as a mode-2 fiber. A slice is a two-dimensional section of a tensor,
  obtained by fixing one index; $\cA(i,:,:)$ is a
  mode-1 slice or 1-slice. A particular element of a 3-tensor $\cA$ can
  be denoted in two different way, i.e., "\matlab-like" notation and
  standard subscripts with $\cA(i,j,k)$ and $a_{ijk}$, respectively.
 \begin{definition}
 A 3-tensor $\cA\in \RR^{m \x m \x n}$ is called (1,2)-symmetric if all its 3-slices are symmetric, i.e.,
 \[\cA(i,j,k)=\cA(j,i,k), \qquad i,j=1,2,\ldots,m, \qquad k=1,2,\ldots,n.\]
 \end{definition}
 Symmetry with respect to any two modes and for tensors of higher
 order than 3 can be defined analogously.
 
 We use $I_k$ for the identity matrix of dimension $k$.

 \subsection{Multilinear tensor-matrix multiplication}
 \label{sec:ttm}
  We  first consider the multiplication of a tensor by a matrix. When
  a tensor is multiplied by a single matrix in mode i, say, we will
  call the 
 operation the mode-i multilinear multiplication (or mode-i product) of a tensor by a
 matrix. For example the mode-1 product of a tensor $\cA\in \RR^{l\x m\x n}$ by a matrix $U\in \RR^{p\x l}$ is defined 
 \[
   \RR^{p\x m\x n}\ni\cB=\tml[1]{U}{\cA},\qquad
   b_{ijk}=\sum_{\a=1}^{l}u_{i\a}a_{\a jk}.
 \]
 This means that all mode-1 fibers in the 3-tensor $\cA$ are
 multiplied by the matrix $U$. The mode-2 and the mode-3 product are
 defined in a similar way.
 Let  the
 matrices $V\in \RR^{q\x m}$ and $W\in \RR^{r\x n}$;  multiplication
 along all three modes is defined
 \[
 \RR^{p\x q\x r}\ni\cB=\tml{U,V,W}{\cA}, \qquad b_{ijk}=\sum_{\a=1}^{l}\sum_{\b=1}^{m}\sum_{\g=1}^{n}
 u_{i\a}v_{j\b}w_{k\g}a_{\a \b\g}.
 \]
 For multiplication with a transposed matrix $X\in\RR^{l\x s}$ it is convenient to introduce a separate notation,
 \[\RR^{s\x m\x n}\ni\cB=\tml[1]{X^T}{\cA}=\tmr[1]{\cA}{X},\qquad
 b_{i,j,k}=\sum_{\a=1}^{l}x_{\a i}a_{\a jk}.\]
 In a similar way if $x\in\RR^l$ then
 \[
   \RR^{1\x m\x
     n}\ni\cB=\tml[1]{x^T}{\cA}=\tmr[1]{\cA}{x}=B\in\RR^{m\x n}.
 \]
Thus the tensor $\cB$ is identified with a matrix $B$.  
 
 \subsection{ Inner product and norm, contractions}
 \label{sec:contraction}
 The \emph{inner product} of two tensors $\cA$ and $\cB$ of the same order and dimensions is denoted by $\<\cA,\cB\>$
 and is computed as a sum of element-wise products over all the indices, that is
 \[
 \<\cA,\cB\>=\sum_{i,j,k}a_{ijk}b_{ijk}
 \]
 The  product allows us to define the  Frobenius norm of a tensor $\cA$ as
 \[
   \norm\cA=\<\cA,\cA\>^{1/2}.
  \]
  As in the matrix case the Frobenius norm of a tensor is invariant
  under orthogonal 
  transformations, i.e., $\norm A=\norm{\tml{U,V,W}{\cA}}$, for orthogonal matrices $U$, $V$, and $W$.
  This follows immediately from the fact that mode-i multiplication by an orthogonal
  matrix does not change the Euclidean length of the mode-i fibers.

 The inner product is a \emph{contraction}. We also define \emph{partial contractions} that involve less than three modes,
 \begin{align*}
   \cC &= \<\cA,\cB\>_1, & c_{jk\mu\nu} &=\sum_{\lambda}a_{\lambda
     jk}b_{\lambda\mu\nu},\\
   \cD&=\<\cA,\cB\>_{1,2}, & d_{k\nu}&=\sum_{\lambda,\mu}a_{\lambda\mu k}b_{\lambda\mu\nu}.
 \end{align*}
 We use negative subscripts to denote partial contractions in all but
 one mode,
 \begin{equation}
   \label{eq:partial-c}
 \<\cA,\cB\>_{-1}= \< \cA, \cB \>_{2,3}.   
 \end{equation}
 The result is a matrix (order 2 tensor) of inner products between the
 mode-1 slices of the two tensors. 
 For partial contractions only the contracted modes are required to be
 equal, so the result matrix may be rectangular.

 \subsection{Multilinear rank}

 Unlike the matrix case, there is no
 unique definition of the rank of a tensor. In this paper we consider
 the concept of  multilinear rank defined by Hitchcock
 \cite{hitchcock1928multiple}. Let $A^{(i)}$ denote the mode-$i$
 unfolding (matricization) of $\cA$ (using some ordering of the vectors),
 \[A^{(i)}=\mathtt{unfold}_i(\cA),\]
  where the columns of $A^{(i)}$ are  all  mode-$i$
  fibers \cite{de2000multilinear}. Similarly, let
      $\mathtt{fold}_i$ be the inverse of $\mathtt{unfold}_i$. The multilinear rank of a third order tensor $\cA$ is an
 integer triplet $(p, q, r)$ such that
\[   p=\mathrm{rank}(A^{(1)}), \qquad 
   q=\mathrm{rank}(A^{(2)}),\qquad 
   r=\mathrm{rank}(A^{(3)}),
 \]
 where  $\mathrm{rank}(A^{(i)})$ is the matrix rank. In this paper we will deal only with
 multilinear rank, and we will use the notation rank-$(p, q, r)$, and
      $\mathrm{rank}(A) = (p, q, r)$. 
 For matrices the rank is obtained via the \svd; see, e.g.,
 \cite[Chapter 2]{golub2013matrix}. In exact 
 arithmetic the multilinear rank can be computed using the higher order singular value
      decomposition (HOSVD) \cite{de2000multilinear}.

 \subsection{Best rank-$(r_1,r_2,r_3)$ approximation}\label{sec:best rank}
 The problem \eqref{eq:min} of approximating a give\-n tensor $\cA\in
      \RR^{l\x m\x n}$ by another tensor $\cB$ of equal dimensions but
      of lower rank, 
 occurs in many modern applications, e.g., machine learning
      \cite{lim2008cumulant}, pattern classification
      \cite{savas2007handwritten}, analytical and quantum chemistry
      \cite{smilde2005multi, khoromskij2007low}, and signal processing
      \cite{comon1996decomposition}. We assume that
      $\mathrm{rank}(\cB)=(r_1,r_2,r_3)$, which means that $\cB$ can be
      written as a product of a core tensor $\cF\in\RR^{r_1\x r_2\x r_3}$
      and three matrices, 
 \[
 \cB=\tml{U,V,W}{\cF},\qquad
 \cB(i,j,k)=\sum_{\l,\mu,\nu=1}^{r_1,r_2,r_3}u_{i\l}v_{j\mu}w_{k\nu}f_{\lambda\mu\nu}, 
 \]
 where $U\in\RR^{l\x r_1}$, $V\in\RR^{m\x r_2}$, and $W\in\RR^{n\x r_3}$ are
 full column rank matrices. Without loss a generality, we can suppose
 that $U$, $V$ and $W$ have orthonormal columns, as any
 nonorthogonality may be incorporated into $\cF$\footnote{Assume
   that $U=U_0 R_0$ is the thin QR decomposition of $U$. Then
   $\tml{U,V,W}{\cF} = \tml{U_0 R_0,V,W}{\cF} =
   \tml{U_0,V,W}{(\tml[1]{R_0}{\cF})} =:  \tml{U_0,V,W}{\cF_0}$.}. Therefore the best
 multilinear low rank problem (\ref{eq:min}) can be written as 
 \begin{equation}\label{eq:main min}
 \min\limits_{\cF,U,V,W}\norm{\cA-\tml{U,V,W}{\cF}},\quad
 \mathrm{subject~  to} \quad  U^TU=I_{r_1},\quad V^TV=I_{r_2},\quad W^TW=I_{r_3}.
 \end{equation}
 There are a few major differences
 between the best low rank approximation of matrices and 
      3-mode tensors  and higher. In the matrix case, the explicit solution of
      corresponding problem can be obtained from the  SVD, see the
      Eckart-Young property in \cite[Theorem 2.4.8]{golub2013matrix}. A simple
      proof is given in \cite[Theorem 6.7]{elden2007matrix}. There is
      no known closed form solution for the minimization problem
      (\ref{eq:main min}), but it can be shown that this is a well-defined
      problem in the sense that for any $(r_1,r_2,r_3)$ a solution
      exists \cite[Corollary 4.5]{de2008tensor}. Several iterative 
      methods for computing the low rank approximation for small and
      medium size tensors have been
      proposed, see \cite{de2000best, elden2009newton,
      ishteva2009differential, savas2010quasi}. In \cite{de2000best}, it is
      shown that (\ref{eq:main
      min}) is equivalent to following maximization problem
 \begin{equation}\label{eq:max}
 \max\limits_{U,V,W} \Phi(U,V,W),\quad \mathrm{subject~  to}\quad
 U^TU=I_{r_1},\quad V^TV=I_{r_2},\quad W^TW=I_{r_3}, 
 \end{equation}
      where $\Phi(U,V,W)=\norm{\tmr{\cA}{U,V,W}}^2$. Since the  norm is
      invariant under orthogonal transformations, it holds that
      $\Phi(U,V,W)=\Phi(UQ_1,VQ_2,WQ_3)$ for any orthogonal matrices
      $Q_1\in\RR^{r_1\x r_1}$, $Q_2\in\RR^{r_2\x r_2}$ and $Q_3\in\RR^{r_3\x
        r_3}$. Hence  \eqref{eq:max} is       equivalent to a
      maximization problem  
       over a product of Grassmann manifold; for optimization on
      matrix manifolds, see
      \cite{edelman1998geometry,ams:07,ishteva2011best,ishteva2009differential}. 

      After computing the optimal $U$, $V$ and $W$ the optimal $\cF$
 can be obtained by considering the minimization of (\ref{eq:main min}) as a
      linear least squares problem with unknown  $\cF$. 
      
 \begin{lemma}
 Let $\cA\in\RR^{l\x \x m\x n}$ be given along with three matrices
 with orthonormal columns, $U\in\RR^{l\x r_1}$, $V\in\RR^{m\x r_2}$, and
 $W\in\RR^{n\x r_3}$, where $r_1\leq l$, $r_2\leq m$, and $r_3\leq n$. Then
 the least squares problem 
 \[\min\limits_{\cF}\norm{\cA-\tml{U,V,W}{\cF}}\]
 has the unique solution
 \begin{equation}
   \label{eq:Rayleigh}
 \cF=\tml{U^T,V^T,W^T}{\cA}=\tmr{\cA}{U,V,W}.   
 \end{equation}
 %
\end{lemma}

       For a proof, see       e.g. \cite{de2000best,savas2013krylov}. 
 The tensor $\cF$ is a generalization of the matrix Rayleigh
      quotient.

 \subsection{Gradient on the product manifold}\label{subsec:gradient}

      In \cite{elden2009newton} a Newton-Grass\-mann me\-thod is derived
      for computing the solution of maximization problem
      \eqref{eq:max}. The constraints on the unknown matrices $U$, $V$,
      and $W$ are taken into 
      account by formulating the problem as an optimization problem on a
      product of three      Grassmann manifolds.  In this paper we
      will need the gradient of $\Phi$ in the tangent space of the product
      manifold for a stopping criterion. 
      This  gradient can be expressed   in the ambient coordinate
      system, or in       local coordinates. 
 In the context of the new methods presented  it is practical and more
 efficient  to        use local coordinate 
 representations, see Proposition \ref{prop:gradient}. Let $(U \, U_\perp)$ denote the enlargement of $U$ to
 a square orthogonal matrix, and use the analogous notation for $V$
 and $W$. Then  the Grassmann
 gradient at $(U,V,W)$  can
 be written as  
 \begin{equation*}
   \label{eq:grad-local}
   \nabla_{\mathrm{local}}(U,V,W) = (\<\cF^1_\perp,\cF \>_{-1},
   \<\cF^2_\perp,\cF \>_{-2},\<\cF^3_\perp,\cF \>_{-3}),
 \end{equation*}
 where $\cF^1_\perp = \tmr{\cA}{U_\perp,V,W}$,  $\cF^2_\perp =
   \tmr{\cA}{U,V_\perp,W}$,   and $\cF^3_\perp =
     \tmr{\cA}{U,V,W_\perp}$.  In the context of the HOOI, see Section
     \ref{sec:HOOI}, it is more 
      efficient to use global 
      coordinates. For instance, the first component of the gradient
      can be computed       as 
      \[
        (I- U U\tp) \Gamma_1=\Gamma_1 - U (U\tp \Gamma_1), \qquad \Gamma_1 =
        \langle \tmr[2,3]{\cA}{V,W} , \cF \rangle_{-1}.
        \]        
 For more details on coordinate
      representations for this problem, see
      \cite{elden2009newton}, \cite[Section
      3.2]{elden2011perturbation}.
In the rest of this paper, the concept \emph{G-gradient} will mean the
Grassmann gradient in global or local coordinates.

\subsection{Conditioning of the best approximation problem}
\label{sec:conditioning}

      The best rank-$r$ approximation problem for a matrix $A$ is not
      unique if the singular values satisfy $\sigma_r(A) =
      \sigma_{r+1}(A)$. The problem is ill-conditioned if the
      \emph{gap} is small, i.e. $\sigma_r(A) > \sigma_{r+1}(A)$ but
      $\sigma_r(A) \approx \sigma_{r+1}(A)$, see e.g. \cite[Chapter
      3]{stewart2001matrix}, \cite[Chapter 8.6]{golub2013matrix}. 
      A similar situation  exists for the tensor case
      \cite[Corollary 4.5]{elden2011perturbation}  (note that the
      perturbation theory for the SVD is a special case of that for
      the best rank-$(r_1,r_2,r_3)$ approximation of a
      tensor).  Define
       \begin{align*}
         s_i^{(k)} &= (\lambda_i(\langle \cF, \cF \rangle_{-k}))^{1/2},
                            \qquad i=1,2,\ldots,r_k, \quad k=1,2,3,\\
         s_{r_{k}+1}^{(k)} &= (\lambda_{\max}(\langle \cF_\perp^k,
        \cF_\perp^k \rangle_{-k}))^{1/2}, \qquad k=1,2,3,
        \end{align*}
        where the $\lambda$'s are eigenvalues, in descending order, of
        the symmetric         matrices.  We will refer to these
        quantities as 
        \emph{S-values}. 
        Then we can define three gaps, one for each mode, 
        \[
          \gap_k = s_{r_k}^{(k)} - s_{r_{k}+1}^{(k)}, \qquad
          k=1,2,3.
        \]
        (In the matrix case,  there is only one set of  $s_i^{(k)}$,
        which  are the singular         values).  It is shown in
        \cite[Section 5.3]{elden2011perturbation} that the gaps can be
        taken as 
      measures of the conditioning of the  best approximation 
      problem.  If, for any $k$, $s_{r_k}^{(k)}$ is considerably larger than
      $s_{r_{k}+1}^{(k)}$  then the approximation problem is
      well-conditioned with respect to  mode $k$. Conversely,
      if the gap is small, then the problem is ill-conditioned.

      \section{Krylov methods for  matrices}\label{sec:krylov-type}

       Krylov
      subspace methods  are the main  class of algorithms
      for solving iteratively large and sparse matrix problems. Here
      we give a brief  introduction to  Krylov methods, illustrating
      with the Arnoldi method. 

For a  given  square matrix $A\in\RR^{n\x n}$ and a nonzero vector
$u\in\RR^n$  the subspace 
 \begin{equation}\label{eq:krylov}
 \mathcal{K}_k(A,u)=\mathrm{span}\{u,Au,A^2u,\ldots,A^{k-1}u\}
 \end{equation}
      is called the  Krylov subspace of dimension $k$ associated with
      $A$ and $u$ \cite[Chapter 10]{golub2013matrix}. The Arnoldi
      process \cite[Chapter 10]{golub2013matrix}
      which is obtained by  applying Gram-Schmidt orthogonalization,
      generates an orthogonal basis for the Krylov subspace
      (\ref{eq:krylov}). The recursive Arnoldi process is equivalent to the
       matrix equation
 \begin{equation}\label{eq:arnoldi}
 AU_k=U_{k+1}\hat{H_k},
 \end{equation}
 where $U_k=[u_1,\ldots,u_k]$ is an orthogonal matrix and $\hat{H_k}\in\RR^{(k+1)\x k}$ is an Hessenberg matrix with orthogonalization coefficients.
  Using this factorization one can compute an approximate solution of a
  linear system or an eigenvalue problem with coefficient matrix $A$ by
  projecting onto the  low dimensional subspace $U_k$. 

 \subsection{The Krylov-Schur approach to computing  low rank
   approximations}\label{sec:K-Schur}

 The Arnoldi and Lanczos methods (see, e.g., \cite[Chapter
 10]{golub2013matrix}, \cite[Chapter 5]{stewart2001matrix})
 are
 designed for computing the eigenvalues of  large and sparse
 matrices. Writing the Krylov decomposition \eqref{eq:arnoldi} in the
 form
 \begin{equation}\label{eq:arnoldi2}
 AU_k=U_kH_k+\b_ku_{k+1}e_k^T,
 \end{equation}
 approximations are obtained by computing the eigenvalues of $H_k$ (the
 Ritz values),
 using an algorithm for dense matrices. A problem with this approach is
 that when $k$ grows the cost for orthogonalizing a newly generated
 vector against the columns of $U_k$ increases. In addition, since
 $U_k$ is dense, the memory requirements may increase too much. In
 order to save memory and work, an implicit restarting Arnoldi technique
 (IRA) was developed \cite{sorensen1992implicit}, and implemented in the
 highly successful ARPACK package \cite{lehoucq1998arpack}.
  In \cite{stewart2002krylov} Stewart proposed a
 Krylov-Schur method which is mathematically equivalent to IRA.
 In the following we give
 a brief introduction to the Krylov-Schur method, based on
 \cite[Chapter 5]{stewart2001matrix}.

 Consider  the Krylov
 decomposition \eqref{eq:arnoldi2} of $A$, which  is obtained after $k$ steps
 of the Arnoldi recursion.
 Let  $k=r+k_2$ where $r$ is the number of $A$'s eigenvalues that
 we want to compute. The Schur decomposition of $H_k$ is
 \begin{equation}
   \label{eq:schur-Hk}
   H_k = Q T_k Q^T,
 \end{equation}
 where $Q$ and $T_k$ are  orthogonal and upper triangular matrices,
 respectively\footnote{In real arithmetic $T_k$ may have $2 \times 2$
   blocks on the diagonal, corresponding to complex conjugate
   eigenvalues. For simplicity of presentation we ignore this
   here.}.  Partition  $T_k$ and $Q$ as 
 \begin{equation}
   \label{eq:Tk-part}
 T_k=\begin{pmatrix}
       T_{11} & T_{12} \\
       0 & T_{22} \\
     \end{pmatrix},
 \qquad Q=[Q_1\;Q_2],  
 \end{equation}
 where $Q_1\in\RR^{k\x r}$, $Q_2\in\RR^{k\x k_2}$, and the
 eigenvalues have been ordered such that those of
 $T_{11}\in\RR^{r\x r}$ are of interest and those of
 $T_{22}\in\RR^{k_2\x k_2}$ are not. The Krylov decomposition
 (\ref{eq:arnoldi2}) can now be written
 \begin{equation}\label{eq:k-s}
 A(U_kQ_1\; U_kQ_2)=(U_kQ_1\; U_kQ_2)\begin{pmatrix}
       T_{11} & T_{12} \\
       0 & T_{22} \\
     \end{pmatrix}+u_{k+1}(b_{1,k+1}^T\;b_{2,k+1}^T),
 \end{equation}
 where $(b_{1,k+1}^T\;b_{2,k+1}^T)=(\b_k e_k^TQ_1\;\b_k e_k^TQ_2)$.
 Then
 \begin{equation}\label{ksd}
 A\hat{U}_{r}=\hat{U}_{r}T_{11}+u_{k+1}b_{1,k+1}^T,
 \end{equation}
      where $\hat{U}_{r}=U_kQ_1$, is also a Krylov decomposition.
      The main idea of the Krylov-Schur algorithm is to repeatedly expand
      (\ref{ksd}) to dimension $k$, using the Arnoldi recursion, and
      then reduce it back to 
      dimension $r$. The dimension of $H_k$ is assumed to be much
      smaller than that of $A$, which makes it possible to compute the
      Schur decomposition of $H_k$ using the dense QR algorithm
      (implemented in LAPACK).
      A sketch of this method \cite[Chapter 5]{stewart2001matrix} is
      given in Algorithm \ref{alg:M1}.

      \begin{algorithm}[htb]
 {\vskip 4pt}
 \begin{algorithmic}
 {\vskip 4pt}
 \STATE {Given: matrix $A$, the number of desired eigenvalues $r$ and the number of steps $k$.}
 \STATE {Build a Krylov decomposition (\ref{eq:arnoldi2}) of $A$.} 
 \STATE {\textbf{Until convergence}}
 \STATE \quad{Compute the Schur decomposition \eqref{eq:schur-Hk} and
   partition as in \eqref{eq:Tk-part} with the $r$}
  \STATE\quad\quad {desired eigenvalues in the block $T_{11}$.} 
 \quad\STATE \quad{Discard the unwanted eigenvalues and save the truncated Krylov-Schur}\STATE\quad\quad{decomposition (\ref{ksd}).}
 \quad\STATE \quad{Check convergence}
 \quad\STATE \quad{Expand  the Krylov decomposition (\ref{ksd}) to one
   of dimension $k$.}
  \end{algorithmic}
  \caption{The Krylov-Schur algorithm}
  \label{alg:M1}
 \end{algorithm}

 The convergence check is based on the eigenvalue residual
 \cite[Chapter 5]{stewart2001matrix}.

 \section{Krylov-type methods for tensors}\label{sec:Krylov-tensor}
      Krylov-type
      methods that generalize the Arnoldi method to tensors are proposed in
      \cite{savas2013krylov}.  The methods are called 
      Krylov-type methods,      because the recursions are
      generalizations of matrix Krylov recursions, but no analogues of
            Krylov subspaces can be identified (to our knowledge).
            These methods are also 
      inspired by  Golub-Kahan bidiagonalization
      \cite{golub1965calculating}.  The
      bidiagonalization process
      generates two sequences of orthonormal basis vectors for certain Krylov
      subspaces. In the tensor case   three sequences of orthogonal
      basis vectors  are generated that are      used to compute a
      core tensor that corresponds to the matrix $H_k$ in
      \eqref{eq:arnoldi2}. For matrix Krylov methods, once an initial
      vector has been
      selected, the following vectors are determined uniquely; in the
      tensor case, one can choose different combinations of  previously
      computed basis vectors. So there
      are different variants of tensor Krylov-type methods.  We
      describe briefly two examples in the following.

 \subsection{Minimal Krylov recursion}\label{sec:min-krylov}
      For a given third order tensor $\cA\in\RR^{l\x m\x n}$ and starting
      two vectors, $u_1\in\RR^l$ and $v_1\in\RR^m$, we can obtain a third
      mode vector by $w_1=\tmr[1,2]{\cA}{u_1,v_1}\in\RR^n$. By using the
      most recently obtained vectors,  three
      sequences of vectors can be generated.
 Using the modified Gram-Schmidt process, a newly generated vector is
 immediately orthogonalized against all the previous ones in its
 mode. The minimal Krylov recursion \cite[Algorithm
      3]{savas2013krylov},which can be seen as a generalization of the
 Golub-Kahan bidiagonalization method,  is given in Algorithm
 \ref{alg:min-krylov}.  The $\mathtt{orth}$ function orthogonalizes $\hat{u}$
      against $U_i$, and normalizes it.
 %
 \begin{algorithm}[htb]
 {\vskip 4pt}
 \begin{algorithmic}
 {\vskip 4pt}
 \STATE{Given: two normalized starting vectors $u_1$ and $v_1$,}
 \STATE{$\a_w w_1=\tmr[1,2]{\cA}{u_1,v_1}$}
 \FOR {$i=1,2,\ldots,k-1$}
     \STATE{$\hat{u}=\tmr[2,3]{\cA}{v_i,w_i}; \quad%
     u_{i+1}=\mathtt{orth}(\hat u,U_i); \quad U_{i+1}=[U_i\;u_{i+1}]$ }
 {\vskip 4pt}
     \STATE{$\hat{v}=\tmr[1,3]{\cA}{u_{i+1},w_i}; \quad%
      v_{i+1}=\mathtt{orth}(\hat v,V_i); \quad V_{i+1}=[V_i\;v_{i+1}]$ }
 {\vskip 4pt}
     \STATE{$\hat{w}=\tmr[1,2]{\cA}{u_{i+1},v_{i+1}}; \quad
        w_{i+1}=\mathtt{orth}(\hat w,W_i); \quad W_{i+1}=[W_i\;w_{i+1}]$ }
 \ENDFOR
   \end{algorithmic}
   \caption{Minimal Krylov recursion }
   \label{alg:min-krylov}
 \end{algorithm}
      Using the three orthogonal matrices $U_k$, $V_k$, and $W_k$
      generated  by Algorithm \ref{alg:min-krylov}, we  obtain a
      rank-$k$  approximation of 
      $\cA$ as
      \[
      \cA\approx\tml{U_k,V_k,W_k}{\cH},\quad\cH=\tmr{\cA}{U_k,V_k,W_k}\in\RR^{k\x
      k\x k}.
      \]

      \subsection{Maximal Krylov recursion}\label{max-krylov}
      In the minimal Krylov recursion a new vector $u_{i+1}$ is generated
      based on the two most recently computed $v_i$ and $w_i$ (and
      correspondingly  for the other modes).  But we can choose any  other
      available $v_j$ and $w_k$ that have not been combined before in an
      operation $\tmr[2,3]{\cA}{v_j,w_k}$.   If we decide to use
      all available combinations, then this is called the Maximal Krylov
      recursion      \cite[Algorithm 5]{savas2013krylov}. Given $V_j$
      and $W_k$, all combinations can be computed as
      \[
      \hat U=\tmr[2,3]{\cA}{V_j,W_k}.
      \]
      Next the mode-1 fibers of the tensor $\hat U$ have to be orthogonalized.
      The number of basis vectors generated      grows very
      quickly. In the following subsections we will describe a few
      modifications  of the maximal recursion that avoid computing many of
      the vectors in the maximal recursion, while maintaining as much
      as possible of the      approximation performance.


      \subsection{Block-Krylov methods}
      \label{sec:block-krylov-methods}
      For large and sparse tensors it is convenient to use software
      that implements operations with tensors, in particular
      tensor-matrix multiplication.  In our numerical experiments  we use
      the extension of 
      the Matlab tensor toolbox \cite{bader2006algorithm} to sparse tensors
      \cite{bako:08}; it is natural to  base some algorithm 
      design decisions  on the use of such software. Other languages
      and implementations are likely to have similar properties. 

      There are two main reasons why we choose to use block-Krylov
      methods. Firstly, it is easier to design and describe
      modifications of the 
      maximal Krylov recursion in terms of blocks. Secondly, and more
      importantly,  the time       required for the 
      computation of sparse tensor-vector and tensor-matrix products is
      dominated by data movements  \cite[Sections
      3.2.4-3.2.5]{bako:08}, where the tensor 
      is reorganized to a different format before the multiplication
      takes place. Therefore, it is better to reuse the reorganized
      tensor for several vectors, as is done in a tensor-matrix
      product (akin to the  use of BLAS 3 operations in dense linear
      algebra). In our experiments with a few sparse tensors of
      moderate       dimensions
      (approximately $500 \times 500 \times 500$ and $3600 \times 3600
      \times 60$) and  6 vectors,  the time for repeated
      tensor-vector multiplication was 3-9 times longer than for the
      corresponding tensor-matrix block multiplication. One should keep in
      mind, however, that  such  timings
      may also be highly problem-, software- and  hardware-dependent.

      Let $\cA\in\RR^{l \x m \x n}$ be given, and  assume that,
      starting from $U_0 \in \RR^{l \x r_1}$, 
      $V_0 \in \RR^{m \x r_2}$ and $W_0 \in \RR^{n \x r_3}$,  three sequences of
      blocks of orthonormal  basis vectors (referred to as ON-blocks) have been 
      computed, $\widehat{U}_{\lambda-1}=[U_0 \,  U_1   \cdots
      U_{\lambda-1}]$, $\widehat{V}_\mu = [V_0
      \, V_1 \, \cdots V_\mu]$,  and $\widehat{W}_{\nu}=[W_0 \, W_1
      \cdots W_{\nu}]$. Letting  $p$ be a block
      size, and 
      $\bar V  \in\RR^{m\x p}$ and $\bar W \in\RR^{n\x p}$ 
      be  blocks selected out of   $\widehat{V}_\mu$ and $\widehat{W}_{\nu}$,
      we compute   a new block $U_\lambda \in\RR^{l\x p^2}$      using
      Algorithm \ref{alg:block-krylov-step}.  
 \begin{algorithm}[htb]
 {\vskip 4pt}
 \begin{algorithmic}
 {\vskip 4pt}
 \STATE{(i). \quad
   ${\cU}^{(1)}=\tmr[2,3]{\cA}{\bar{V},\bar{W}} \in \RR^{l
     \x  p \x p}$,}
 {\vskip 4pt}
     \STATE{(ii). \quad $\cH^{1}_{\lambda-1} = \tml[1]{\widehat{U}_{\lambda-1}^T}{\cU^{(1)}} =
                    \tmr{\cA}{\widehat{U}_{\lambda-1},\bar V, \bar W},$}
 {\vskip 4pt}
     \STATE{(iii). \quad $\widetilde{\cU}^{(1)} = {\cU}^{(1)} -
                                \tml[1]{\widehat{U}_{\lambda-1}}{\cH_{\lambda-1}^{1}}$,}
 {\vskip 4pt}
     \STATE{(iv). \quad $\tml[1]{U_\lambda}{\cH^1_\lambda} = \widetilde{\cU}^{(1)}.$
         \qquad Orthonormalize the mode-1 fibers of $\widetilde{\cU}^{(1)}$.  }
   \end{algorithmic}
   \caption{Generic mode-1 block-Krylov step }
   \label{alg:block-krylov-step}
 \end{algorithm}

      The algorithm is written in tensor form, in order to make the
      operations in Steps (ii)-(iv) look like the Gram-Schmidt
      orthogonalization that it is. In our actual implementations we have
      avoided some tensor-matrix restructurings, see Appendix
      \ref{sec:block-krylov}.  

 Clearly, after Step (iii) we have
      $\tml[1]{\widehat{U}_{\lambda-1}^T}{\widetilde{\cU}^{(1)}} =0$.
      In Step (iv) the mode-1 vectors (fibers) of
 $\widetilde{\cU}^{(1)}$   are orthogonalized and organized in a
 matrix $U_\lambda \in \RR^{l \x p^2}$, i.e. a ``thin'' QR
 decomposition  is computed,
 \[
 \mathtt{unfold_1}( \widetilde{\cU}^{(1)}) = U_\lambda H_\lambda^1,
 \]
  where the matrix  $H_\lambda^1$ is upper triangular.  The
      tensor $\cH^1_\lambda 
 = \mathtt{fold}_1(H_\lambda^1)  \in \RR^{p^2 \x p
      \x p}$ contains the 
   columns of $H_\lambda^1$, and consequently it has a
     ``quasi-triangular'' structure.   

 The mode-1 step can be written
 \begin{equation}
   \label{eq:Mode-1-BK-recursion}
 \tmr[2,3]{\cA}{\bar{V},\bar{W}} =
 \tml[1]{\widehat{U}_{\lambda-1}}{\cH^{1}_{\lambda-1}} +  \tml[1]{U_\lambda}{\cH^1_\lambda}.   
 \end{equation}
 This can be seen as a  generalization to blocks and tensors of one
 step of the Arnoldi recursion, cf.       \eqref{eq:arnoldi2}. 

 The mode-2 and mode-3 block-Krylov steps are analogous. Different
 variants of BK methods can be derived using different choices of
 $\bar V$ and $\bar W$, etc. However, we will always assume that the
 blocks $U_0$, $V_0$, and $W_0$ are used to generate the blocks with
 subscript 1. 

 The recursion \eqref{eq:Mode-1-BK-recursion} and its mode-2 and mode-3
 counterparts imply the following simple lemma\footnote{More general
   results can 
   be obtained for    $\tmr{\cA}{U_\lambda,V_\mu,W_\nu}$. Taken
      together, those results can be used to show the existence of a
      tensor structure that is analogous to the block Hessenberg structure
      obtained in a block-Arnoldi method       for     a matrix.}. It
    will   be useful in the computation of gradients.

 \begin{lemma}\label{lem:rayleigh-struct}
 Assume that, starting from three ON-blocks $U_0$, $V_0$, and $W_0$,
 Algorithm \ref{alg:block-krylov-step} 
  and its mode-2 and mode-3
 counterparts, have been used to generate ON-blocks
 $\widehat{U}_{\lambda}=[U_0 \,  U_1   \cdots 
      U_{\lambda}]$, $\widehat{V}_\mu = [V_0
      \, V_1 \, \cdots V_\mu]$,  and $\widehat{W}_{\nu}=[W_0 \, W_1
      \cdots W_{\nu}]$. Then
      \begin{equation}
        \label{eq:AUjV0W0}
        \tmr{\cA}{U_j,V_0,W_0} = %
        \begin{cases}
          \cH_0^1 =  \tmr{\cA}{U_0,V_0,W_0}, & j=0, \\
          \cH_1^1 = \tmr{\cA}{U_1,V_0,W_0}, & j=1, \\
            0, & 2 \leq j \leq \lambda.
        \end{cases}
      \end{equation}
 The corresponding identities hold for modes 2 and 3.
 \end{lemma}

 \begin{proof}
 Consider the identity
   \eqref{eq:Mode-1-BK-recursion} with 
   $\bar V = V_0$, and $\bar W =W_0$, multiply by $U_j^T$ in the first
   mode, and use orthogonality and
   $\tml[1]{U_j^T}{(\tmr[2,3]{\cA}{V_0,W_0})}=\tmr{\cA}{U_j,V_0,W_0}$. 
 \end{proof}

 \begin{proposition}\label{prop:gradient}
   Let $(U_0,V_0,W_0)$ with $U_0 \in \RR^{l \x r_1}$, $V_0 \in \RR^{m
     \x r_2}$, and $W_0 \in \RR^{n \x r_3}$,   have  orthonormal
   columns.  Let it be a starting point
   for one block-Krylov step in each mode  with $\cA \in \RR^{l \x m
     \x n} $, giving
   $(U_1,V_1,W_1)$ and    tensors
   \begin{align*}
     \cH_0&= \tmr{\cA}{U_0,V_0,W_0}, \qquad 
     \cH_1^1 =\tmr{\cA}{U_1,V_0,W_0},\\
     \cH_1^2&=\tmr{\cA}{U_0,V_1,W_0}, \qquad
          \cH_1^3=\tmr{\cA}{U_0,V_0,W_1}.            
   \end{align*}
   Then the norm of the G-gradient of \eqref{eq:max} at $(U_0,V_0,W_0)$ is
   \[
     \| \nabla(U_0,V_0,W_0) \|^2 =
     \| \langle \cH_0, \cH_1^1\rangle_{-1} \|^2 +
     \| \langle \cH_0, \cH_1^2\rangle_{-2} \|^2 +
     \| \langle \cH_0, \cH_1^3\rangle_{-3} \|^2.
   \]
 \end{proposition}

 \begin{proof}
   The mode-1 gradient at $(U_0,V_0,W_0)$ is $\langle \cF, \cF_\perp^1
   \rangle_{-1}$, where $\cF = \cH_0 = {\tmr{\cA}{U_0,V_0,W_0}}$, and
   $\cF_\perp^1 = \tmr{\cA}{U_\perp,V_0,W_0}$, 
   and $(U_0 \, U_\perp)$ is an orthogonal matrix. So we can write
   $U_\perp = (U_1 \, U_{1\perp})$, where $U_{1\perp}\tp (U_0 \, U_1) =
   0.$  Since from  Lemma  \ref{lem:rayleigh-struct}
   \[
    \RR^{(l-r_1) \x r_2 \x r_3} \ni  \cF_\perp^1 = \tmr{\cA}{(U_1 \, U_{1\perp}),V_0,W_0} =
     \begin{pmatrix}
       \cH_1^1 \\
       0
     \end{pmatrix},
   \]
   the mode-1 result follows. 
   The proof for the other modes is
   analogous. 
 \end{proof}

 Assume we have an algorithm based on block-Krylov steps in all three
 modes, and we want to compute the G-gradient to check if a point
 $(U_0,V_0,W_0)$  is approximately stationary. Then, by Proposition
 \ref{prop:gradient}, we need only perform one block-Krylov step in
 each mode, starting from $(U_0,V_0,W_0)$,   thereby avoiding the
 computation of $U_\perp$, $V_\perp$,  and $W_\perp$, which are
 usually large and dense.  If the norm of the  G-gradient is not small
 enough, then   one would perform more block-Krylov steps. Thus the
 gradient computation comes for  free, essentially. 
 
  \subsubsection{Block-Krylov methods for (1,2)-symmetric tensors}
      \label{sec:block-krylov-sym}

      In many real applications, the tensors are
      $(1,2)$-symmetric. This is the case, for instance, if the
      mode-3 slices      of the tensor represent undirected graphs. In the
      rest of this section we will assume that $\cA \in \RR^{m \x m \x n}$  is
      $(1,2)$-symmetric,  and  we will compute two sequences of blocks
      $U_1,U_2,\ldots$ and $W_1,W_2,\ldots$, where the $U_\lambda$
      blocks  contain  basis vectors for modes 1 and 2, and the $W_\nu$ for
      mode 3.

      A new block $U_\lambda$ is computed from given blocks $\bar U$ and
      $\bar W$ in       the  same way as in the nonsymmetric case,
      using Algorithm       \ref{alg:block-krylov-step}.%

      To compute a  new block $W_\nu$,  we use two blocks $\bar U$ and
      $\bar{\bar U}$. If  $\bar{U} \neq \bar{\bar U}$, then  we can use the
      analogue of       Algorithm  \ref{alg:block-krylov-step}. 
 In the case $\bar{U} = \bar{\bar U}$ the product tensor
      $\tmr[1,2]{\cA}{\bar{U},\bar{U}}$   is
      (1,2)-symmetric, which means that almost half its 3-fibers are
      superfluous, and should be removed. Thus, letting
      $\widetilde{\cW}^{(3)}$ denote the tensor that is obtained in
      (iii) of Algorithm \ref{alg:block-krylov-step},  
      we compute the ``thin'' QR decomposition,
 \[
 \mathtt{triunfold_3}(\widetilde{\cW}^{(3)}) = W_\nu H_\nu^3,
 \]
      where $\mathtt{triunfold_3}(\cX)$ denotes the mode-$3$ unfolding
      of the (1,2)-upper triangular part of the tensor $\cX$. If
      $\bar U \in \RR^{m \x p}$, then $W_\nu \in \RR^{n \x p_\nu}$,
      where $p_\nu=p(p +1)/2$.

      A careful examination of Lemma \ref{lem:rayleigh-struct} for the
      case of (1,2)-symmetric tensors shows that the
      corresponding result holds also here.  We omit the
      derivations  in  order not to make the paper too long. 
                 
      \subsubsection{Min-BK method}\label{sec:min-BK}
      Our  simplest block-Krylov method is the (1,2)-sym\-metric 
 block version of the minimal Krylov recursion of Section
 \ref{sec:min-krylov}, which we 
 refer to as the min-BK method. Here, instead of using only two
 vectors in the multiplication $\hat u = \tmr{\cA}{u_i,w_i}$ we  use
 the $p$       first vectors from the previous  blocks. Let $\bar
 U = U(:,1:p)$, denote the first $p$ vectors of a matrix block $U$.  The
 parameter $s$ in Algorithm \ref{alg:min-BK} is the number of stages. 
 \begin{algorithm}[htb]
 {\vskip 4pt}
 \begin{algorithmic}
 {\vskip 4pt}
 \STATE{\textbf{for }$i=1:s$}
 {\vskip 4pt}
 \STATE{ $\qquad %
   \cU^{(i)} = \tmr[2,3]{\cA}{\bar U_{i-1}, \bar W_{i-1}}$}               
             {\vskip 4pt}
\STATE{$\qquad\qquad \cdots \qquad \% \mbox{ Gram-Schmidt giving}\quad U_{i}$}
{\vskip 4pt}
\STATE{$\qquad \cW^{(i)} = \tmr[1,2]{\cA}{\bar U_{i-1}, \bar U_{i-1}}$}
{\vskip 4pt}
\STATE{$\qquad\qquad \cdots \qquad \% \mbox{ Gram-Schmidt giving}\quad
  W_{i}$ }
{\vskip 4pt}
\STATE{\textbf{end for}}
   \end{algorithmic}
   \caption{Min-BK method for (1,2)-symmetric tensor $\cA$. Given $U_0$  and $W_0$ }
   \label{alg:min-BK}
 \end{algorithm}

 \begin{table}[h]
 \begin{center}
 \caption{\label{table:min-BK} Diagrams of the min-BK method. Top:
   combinations of blocks  for computing the new blocks. The stages
   are separated by solid lines. 
   Bottom:   The number of basis vectors in the stages 
    when $U_0\in\RR^{m\x 2}$,
      $W_0\in\RR^{n\x 2}$, $U_0\in\RR^{m\x 7}$,
      $W_0\in\RR^{n\x 7}$, and $U_0\in\RR^{m\x 10}$, $W_0\in\RR^{n\x
      10}$, respectively, with $p=4$ (separated by $/$).  }

 \end{center}
 \begin{center}
 \begin{tabular}{c|c|c|c|c|}
 \multicolumn{1}{r}{}&\multicolumn{1}{r}{$W_0$}& \multicolumn{1}{r}{$\bar{W}_1$} & \multicolumn{1}{r}{$\bar{W}_2$}&\multicolumn{1}{r}{\ldots}\\
 \cline{2-5}
 $U_0$& $U_1$ & &&\\
 \cline{2-2}
 $\bar{U}_1$& \multicolumn{1}{r}{} & $U_2$&&\\
 \cline{2-3}
 $\bar{U}_2$& \multicolumn{1}{r}{} & \multicolumn{1}{r}{}&$U_{3}$&\\
 \cline{2-4}
 $\vdots$& \multicolumn{1}{r}{} & \multicolumn{1}{r}{}&\multicolumn{1}{r}{}&$\ddots$\\
 \cline{2-5}
 \end{tabular}
 \quad
 \begin{tabular}{c|c|c|c|c|}
 \multicolumn{1}{r}{}&\multicolumn{1}{r}{$U_0$}& \multicolumn{1}{r}{$\bar{U}_1$} & \multicolumn{1}{r}{$\bar{U}_2$}&\multicolumn{1}{r}{\ldots}\\
 \cline{2-5}
 $U_0$& $W_1$ & &&\\
 \cline{2-2}
 $\bar{U}_1$& \multicolumn{1}{r}{} & $W_2$&&\\
 \cline{2-3}
 $\bar{U}_2$& \multicolumn{1}{r}{} & \multicolumn{1}{r}{}&$W_{3}$&\\
 \cline{2-4}
 $\vdots$& \multicolumn{1}{r}{} & \multicolumn{1}{r}{}&\multicolumn{1}{r}{}&$\ddots$\\
 \cline{2-5}
 \end{tabular}
 \end{center}
 \bigskip
 \begin{center}
 \begin{tabular}{|c||c|c||c|c|}
 \hline
 Stage&  & $k_1$&&$k_3$\\
 \hline
 1& $[U_0,U_1]$& 6/56/110&$[W_0,W_1]$&5/35/65\\
 2& $[U_0,\ldots,U_2]$& 18/72/126&$[W_0,\ldots,W_2]$&15/45/75\\
 3& $[U_0,\ldots,U_3]$& 34/88/142&$[W_0,\ldots,W_3]$&25/55/85\\
 \hline
 \end{tabular}
 \end{center}
 \end{table}

 The min-BK method is further described  in the  three diagrams of Table
 \ref{table:min-BK}. Note that to conform with Proposition
 \ref{prop:gradient} we always use $\bar U_0 = U_0$ and $\bar W_0
 =W_0$.  It is seen 
 that the growth of the number of 
 basis vectors, the $k_i$  parameters,  is relatively slow. However, it 
 will be seen in Section \ref{sec:numerical test} that  this method
 is not competitive. 

      \subsubsection{Max-BK method}\label{sec:max-bk}
      The max-BK  method  is maximal in the sense
      that in each stage we use all the available blocks to
      compute new blocks. The algorithm is defined by three diagrams, see Table
      \ref{table:max-bk}. E.g.,  in  stage 2 we use
       combinations of the whole blocks $U_0$, $U_1$, $W_0$, and
       $W_1$, to compute $U_2$, $U_3$, and $U_4$ ($U_1$ was computed
       already in stage 1). 
\begin{table}[h]
 \caption{\label{table:max-bk} Diagrams of the max-BK method. Top: Combinations of blocks of basis vectors for computing  new blocks.
   Bottom:  Number of  basis vectors in the stages of the max-BK
   method. Columns 3 and 5 give the number of basis vectors for
   $r=(2,2,2)$ and $r=(7,7,7)$.}
 \begin{center}
 \begin{tabular}{c|c|c|cc|}
 \multicolumn{1}{r}{}&\multicolumn{1}{r}{$W_0$}& \multicolumn{1}{r}{$W_1$} & \multicolumn{1}{r}{$W_2$}&\multicolumn{1}{r}{$W_3$}\\
 \cline{2-5}
 $U_0$& $U_1$ & $U_2$&$U_5$&$U_{10}$\\
 \cline{2-2}
 $U_1$& \multicolumn{1}{r}{$U_3$} & $U_4$&$U_7$&$U_{12}$\\
 \cline{2-3}
 $U_2$& \multicolumn{1}{r}{$U_6$} & \multicolumn{1}{r}{$U_{8}$}&$U_{9}$&$U_{14}$\\
 $U_3$& \multicolumn{1}{r}{$U_{11}$} & \multicolumn{1}{r}{$U_{13}$}&\multicolumn{1}{r}{$U_{15}$}&$U_{16}$\\
 $U_4$& \multicolumn{1}{r}{$U_{17}$ }& \multicolumn{1}{r}{$U_{18}$}&\multicolumn{1}{r}{$U_{19}$}&$U_{20}$\\
 \cline{2-5}
 \end{tabular}
 \quad
 \begin{tabular}{c|c|c|ccc|}
 \multicolumn{1}{r}{}&\multicolumn{1}{r}{$U_0$}& \multicolumn{1}{r}{$U_1$} & \multicolumn{1}{r}{$U_2$}&\multicolumn{1}{r}{$U_3$}&\multicolumn{1}{r}{$U_4$}\\
 \cline{2-6}
 $U_0$& $W_1$ & $W_2$&$W_4$&$W_7$&$W_{11}$\\
 \cline{2-2}
 $U_1$& \multicolumn{1}{r}{} & $W_3$&$W_5$&$W_8$&$W_{12}$\\
 \cline{2-3}
 $U_2$& \multicolumn{1}{r}{} & \multicolumn{1}{r}{}&$W_{6}$&$W_{9}$&$W_{13}$\\
 $U_3$& \multicolumn{1}{r}{} & \multicolumn{1}{r}{}&\multicolumn{1}{r}{}&$W_{10}$&$W_{14}$\\
 $U_4$& \multicolumn{1}{r}{}& \multicolumn{1}{r}{}&\multicolumn{1}{r}{}&\multicolumn{1}{r}{}&$W_{15}$\\
 \cline{2-6}
 \end{tabular}
 \end{center}
 \bigskip
 \begin{center}
 \begin{tabular}{|c||c|c||c|c|}
 \hline
 Stage&  & $k_1$&&$k_3$\\
 \hline
 1& $[U_0,U_1]$& 6/56&$[W_0,W_1]$&5/35\\
 2& $[U_0,\ldots,U_4]$& 32/1967&$[W_0,\ldots,W_3]$&23/1603\\
 3& $[U_0,\ldots,U_{20}]$& 738/3153108&$[W_0,\ldots,W_{15}]$&530/1935535\\
 \hline
 \end{tabular}
 \end{center}
 \end{table}
      The diagram for the $W_i$'s is triangular
      due to  the (1,2)-symmetry of $\cA$: $\tmr[1,2]{\cA}{U_0,
        U_1}=\tmr[1,2]{\cA}{ U_1,U_0}$, for instance.

      It is clear that the fast  growth of the number of basis vectors
      makes this variant impractical, except for small values of $r$
      and $s$,     e.g. $r=(2,2,2)$ and $s=2$.  In 
      the same way as in the matrix Krylov-Schur method,  we are not
      interested in too
      large values  of $k_1$ and $k_3$, because we will project the
      original problem to one of dimension $(k_1,k_1,k_3)$, which will
      be solved by methods for dense tensors.  Hence we will in the
      next subsection 
      introduce a compromise between the min-BK and  the max-BK method.

       \subsubsection{BK method}\label{sec:mod-bk}
       The BK method is similar to  min-BK  in that it uses
       only the $p$ first vectors of each block in the block-Krylov
       step. In each stage more new blocks are computed than in
       min-BK, but not as many as in max-BK. Both these features  are
      based on numerical tests, where we investigated the performance
      of the BK method in
      the block Krylov-Schur method to be described in Section
      \ref{sec:EK-Schur-gen}. We found that if we omitted the
      diagonal blocks in the diagrams in  Table \ref{table:max-bk},
      then the convergence of the Krylov-Schur method was only
      marginally impeded.  
      The BK method   is described in the three diagrams of
      Table \ref{table:mod-bk}.

       \begin{table}[h]
 \begin{center}
 \caption{\label{table:mod-bk} Diagrams of the BK method. Top:
   Combinations of blocks of basis vectors for computing the new
   blocks. In the case when $U_0\in \RR^{m \x2}$
   and $p \geq 3$, we let $\bar W_1 = W_1$.  
   Bottom:  Basis blocks in the stages of the BK
   method, the number of vectors for $p=4$, and $r=(2,2,2)$,
   $r=(7,7,7)$, and $r=(10,10,10)$.   }
 \end{center}
 \begin{center}
 \begin{tabular}{c|c|c|c|c|}
 \multicolumn{1}{r}{}&\multicolumn{1}{r}{$W_0$}&
                                                 \multicolumn{1}{r}{${\bar W}_1$} & \multicolumn{1}{r}{$\bar{W}_2$}&\multicolumn{1}{r}{$\bar{W}_3$}\\
 \cline{2-5}
 $U_0$& ${U}_1$ & $U_2$&$U_4$&$U_8$ \\
 \cline{2-2}
 $\bar{U}_1$& \multicolumn{1}{r}{$U_3$} & &$U_6$&$U_{10}$ \\
 \cline{2-3}
 $\bar{U}_2$& \multicolumn{1}{r}{$U_5$} &
 \multicolumn{1}{r}{$U_{7}$}&&$U_{12}$ \\
 \cline{2-4}
 $\bar{U}_3$& \multicolumn{1}{r}{$U_9$} & \multicolumn{1}{r}{$U_{11}$}&\multicolumn{1}{r}{$U_{13}$}&\\
 \cline{2-5}
 \end{tabular}
 \quad
 \begin{tabular}{c|c|c|c|c|}
 \multicolumn{1}{r}{}&\multicolumn{1}{r}{$U_0$}& \multicolumn{1}{r}{$\bar{U}_1$} & \multicolumn{1}{r}{$\bar{U}_2$}&\multicolumn{1}{r}{$\bar{U}_3$}\\
 \cline{2-5}
 $U_0$& $W_1$ & $W_2$&$W_3$&$W_5$\\
 \cline{2-2}
 $\bar{U}_1$& \multicolumn{1}{r}{} & &$W_4$&$W_6$\\
 \cline{2-3}
 $\bar{U}_2$& \multicolumn{1}{r}{} & \multicolumn{1}{r}{}&&$W_7$\\
 \cline{2-4}
 $\bar{U}_3$& \multicolumn{1}{r}{} & \multicolumn{1}{r}{}&\multicolumn{1}{r}{}&\\
 \cline{2-5}
 \end{tabular}
 \end{center}
 \bigskip
 \begin{center}
 \begin{tabular}{|c||c|c||c|c|}
 \hline
 Stage&  & $k_1$&&$k_3$\\
 \hline
 1& $[U_0,U_1]$& 6/56/110&$[W_0,W_1]$&5/35/65\\
 2& $[U_0,\ldots,U_3]$& 20/112/190&$[W_0,\ldots,W_2]$&13/63/105\\
 $3$& $[U_0,\ldots,U_7]$& 64/200/302&$[W_0,\ldots,W_4]$&37/107/161\\
 4& $[U_0,\ldots,U_{13}]$& 100/320/448 &$[W_0,\ldots,W_{7}]$&77/167/233\\
 \hline
 \end{tabular}
 \end{center}
\end{table}


It may happen that one of the dimensions of the tensor is
considerably smaller than the other. Assume that $m \gg n$. Then
after a few stages the number of vectors in the third mode may be
equal to $n$, and no more can be generated in that mode. Then
the procedure is modified in an obvious way (the right diagram
is stopped being used)  so that only vectors in the other modes
(U blocks) are 
generated.       The min-BK and max-BK methods can be modified
analogously. 


 \section{A  Tensor Krylov-Schur like method}\label{sec:EK-Schur-gen} 
 When 
 tensor Krylov-type methods are  used for the computation of low
 multilinear rank  approximations of large and sparse tensors,
 they suffer from the  same weakness as Krylov methods 
 for matrices: the
 computational burden for orthogonalizing the vectors as well as
 memory requirements may become prohibitive.
 Therefore a restarting procedure should be tried. We
 will now describe a  generalization of  the
 matrix Krylov-Schur approach  to a corresponding  tensor
 method. Here we assume that the tensor is non-symmetric.  

 Let $\cA\in\RR^{l\x m\x n}$ be  given third order tensor, for  which we
      want to compute the best rank-$(r_1,r_2,r_3)$
      approximation.   For reference we restate the
      maximization problem,
      \begin{equation}
        \label{eq:maxrA}
        \begin{array}{ll}
         &  U \in \RR^{l \x r_1}, \quad U^T U = I_{r_1}, \\
        \max_{U,V,W}  \| \tmr{\cA}{U,V,W} \|, & V \in \RR^{m \x r_2},
                                                \quad V^T V = I_{r_2},
          \\ 
                     &  W \in \RR^{n \x r_3}, \quad W^T W= I_{r_3}.
        \end{array}
      \end{equation}
      %
        
      Let  $k_1$, $k_2$, and $k_3$ be  integers such that
      \begin{equation}
        \label{eq:rk-assumptions}
        r_2 r_3 \leq k_1\ll  l, \quad r_1 r_3 \leq k_2\ll m, \quad
        r_1 r_2 \leq k_3\ll n, 
      \end{equation}
      and assume that we have computed, using   a BK method,  a
      rank-$(k_1,k_2,k_3)$ approximation
 \begin{equation}\label{eq:tkd}
 \cA\approx\tml{X,Y,Z}{\cC}, \qquad \cC = \tmr{\cA}{X,Y,Z},
 \end{equation}
      where $X\in\RR^{l\x k_1}$, $Y\in\RR^{m\x k_2}$, and $Z\in\RR^{n\x
      k_3}$ are  matrices with orthonormal columns, and $\cC\in\RR^{k_1\x k_2\x
      k_3}$ is a core tensor. 
      With the assumption
      \eqref{eq:rk-assumptions}, we can use an algorithm for dense
      tensors, e.g. a Newton-Grassmann  method
                       \cite{elden2009newton,ishteva2009differential},
                       to solve the projected  maximization problem
      \begin{equation}
        \label{eq:maxrC}
        \begin{array}{ll}
         &\hat U \in \RR^{k_1 \x r_1}, \quad \hat U^T \hat U = I_{r_1}, \\
     \max_{\hat U,\hat V,\hat W}  \| \tmr{\cC}{\hat U,\hat V,\hat{W}} \|, &
                 \hat V \in \RR^{k_2 \x r_2}, \quad 
           \hat V^T \hat V = I_{r_2}, \\
          &\hat W \in \RR^{k_3 \x r_3}, \quad \hat W^T \hat W= I_{r_3}.
        \end{array}
      \end{equation}
      %
%
       From the solution of \eqref{eq:maxrC} we have the     best
            rank-$(r_1,r_2,r_3)$ approximation of  $\cC$, 
 \begin{equation}\label{eq:ng}
 \cC\approx\tml{\hat U, \hat V, \hat W}{\cF},
 \end{equation}
      where $\cF \in \RR^{r_1\x r_2\x r_3}$ is the core tensor. This step is
      analogous to computing and truncating  the Schur decomposition
      of the matrix $H_k$ 
      in the matrix case in  Section \ref{sec:K-Schur}. Combining
      (\ref{eq:ng}) and  (\ref{eq:tkd}) we can write
 \begin{equation}\label{eq:tks1}
 \cA\approx\tml{X,Y,Z}\left({\tml{\hat U, \hat V, \hat W}{\cF}}\right)
  = \tml{U,V,W}{\cF},
 \end{equation}
      where $U =X \hat U\in\RR^{l\x r_1}$, $V=Y\hat V\in\RR^{m\x r_2}$, and
      $W=Z \hat W\in\RR^{n\x r_3}$, with orthonormal columns. Thus
      \eqref{eq:tks1} is a rank-$(r_1,r_2,r_3)$ approximation of $\cA$.
       Then, starting with $U_0=U$, $V_0 = V$,
      and $W_0 = W$,  we can again  expand (\ref{eq:tks1}), using a BK
      method,  to       a rank-$(k_1,k_2,k_3)$ approximation
      \eqref{eq:tkd}.
       A sketch of the tensor Krylov-Schur method is given in
      Algorithm \ref{alg:tks}.

      \begin{algorithm}[htb]
 {\vskip 4pt}
 \begin{algorithmic} 
 {\vskip 4pt}
 \STATE {Given: tensor $\cA$, the triplets $(k_1,k_2,k_3)$ and $(r_1,r_2,r_3)$.}
 \STATE {Compute a rank-$(k_1,k_2,k_3)$  approximation \eqref{eq:tkd}
          using a BK method.}
 \STATE {\textbf{Until convergence}}
 \STATE \quad{(i). Compute the best rank-$(r_1,r_2,r_3)$ approximation
   \eqref{eq:ng} of the core }
 \STATE\quad{\quad\;\; tensor $\cC$, using, e.g.,  the Newton-Grassmann
   algorithm. }
 \quad\STATE \quad{(ii). Build the approximation \eqref{eq:tks1}.}
 \quad\STATE \quad{(iii). Perform one step of a BK method starting
   from $(U_0,V_0,W_0)$. } 
\quad\STATE \quad{(iv).   Compute the G-gradient
         at  $(U_0, V_0, W_0)$ (Proposition \ref{prop:gradient}),  and}
 \STATE\quad{\quad\quad\,     check convergence.} 
 \quad\STATE \quad{(v). Expand  to a  rank-$(k_1,k_2,k_3)$
   approximation \eqref{eq:tkd}, using a BK method.}
  \end{algorithmic}
  \caption{Tensor Krylov-Schur  algorithm}
  \label{alg:tks}
 \end{algorithm}

     The algorithm is an outer-inner iteration. In the outer iteration
      \eqref{eq:maxrA} is projected to 
     the problem \eqref{eq:maxrC} 
     using the bases $(X,Y,Z)$.
%
       Step (i) of Algorithm \ref{alg:tks} is the inner iteration, where
      we solve \eqref{eq:maxrC} using an algorithm for a small, dense
      tensor, e.g., the Newton       algorithm on the Grassmann manifold
      \cite{elden2009newton,ishteva2009differential}. 
       The Newton-Grassmann
      method is analogous to and       has the same basic properties
      as the standard Newton method in       a Euclidean space. 

      Notice that if $\cA$ is a $(1,2)$-symmetric tensor, then all
      aforementioned relations are satisfied with $U=V$ and $X=Y$. So
      equation (\ref{eq:tks1}) transfers to
      \[
      \cA\approx\tml{{U},{U},{W}}{\cF}.
      \]

      \subsection{Convergence of the tensor Krylov-Schur algorithm} 
      \label{sec:conv-KS}

      In the discussion below we
      will make the following  assumption:
      \begin{equation}
        \label{eq:stat-assump}
      \begin{array}{ll}
                      & \textrm{the G-Hessian for \eqref{eq:maxrA} at
               the stationary point is negative }\\
             & \textrm{definite, i.e., the stationary point is a
               strict local maximum.}
      \end{array}
      \end{equation}
       Note that 
               \eqref{eq:stat-assump} is equivalent to assuming that
               the objective function is strictly  convex in a neighborhood of
               the stationary point.

\begin{figure}[H]
   \centering
 \includegraphics[height=10cm,width=10cm]{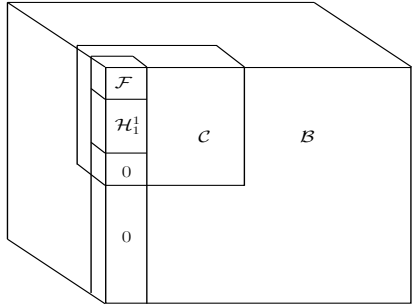}
 \vspace{-4cm}  
 \caption{The tensors $\cB$ and $\cC$, and the  mode-1 subblocks (for
   visibility we have not drawn the corresponding mode-2 and mode-3
   blocks).  The block $\cH_1^1$ is the tensor in the first
   block-Krylov step, see Lemma \ref{lem:rayleigh-struct}.  } 
   \label{fig:CB}
 \end{figure}

 Let $(U_0,V_0,W_0)$ be an approximate solution, and let the
 expanded bases of ON-blocks  be  $ X=[U_0 \, U_1]  \in \RR^{l \x
   k_1}$,        $Y=[V_0\, V_1]  \in \RR^{m \x k_2}$, and $Z=[W_0 \,
 W_1] \in  \RR^{n \x k_3}$. For simplicity we here have included 
   more than one block-Krylov steps in $U_1$, $V_1$, and $W_1$.  Let
   $X_\perp$  be a matrix such that $( X \;  X_\perp)$ 
 is orthogonal, and make the analogous definitions for $Y_\perp$ and
 $Z_\perp$.   We then make the change of  variables 
      \[
        \cB = \tmr{\cA}{[ X \, X_\perp], [ Y \, \
          Y_\perp], [Z \, Z_\perp]}. 
      \]
     The tensor
      \[
        \cC = \tmr{\cA}{ X ,  Y ,  Z }, 
      \]
      is a subtensor of $\cB$, see Figure \ref{fig:CB}. 

      In the discussion of convergence we can, without loss of
      generality, consider the equivalent maximization problem for
      $\cB$,
      \begin{equation}
        \label{eq:maxB}
        \max \| \tmr{\cB}{U,V,W} \|, \qquad U\tp U = I_{r_1}, \quad
        V\tp V = I_{r_2}, \quad W\tp W = I_{r_3}. 
      \end{equation}
      Now  the approximate solution
      $(U_0,V_0,W_0)$ is  represented by
       \[
    E_0 = \left(
     \begin{pmatrix}
       I_{r_1} \\ 0
     \end{pmatrix},
     \begin{pmatrix}
       I_{r_2} \\ 0
     \end{pmatrix},
     \begin{pmatrix}
       I_{r_3} \\ 0
     \end{pmatrix}
   \right), 
 \]
and we have enlarged $E_0$ by one or more block-Krylov steps to 
\[
     E = \left(
     \begin{pmatrix}
       I_{k_1} \\ 0
     \end{pmatrix},
     \begin{pmatrix}
       I_{k_2} \\ 0
     \end{pmatrix},
     \begin{pmatrix}
       I_{k_3} \\ 0
     \end{pmatrix}
   \right).  
 \]
In the inner iteration we shall now compute the best
rank-$(r_1,r_2,r_3)$ approximation for $\cC$,
      \begin{equation}
        \label{eq:maxC}
        \max \| \tmr{\cC}{P,Q,S} \|, \qquad P\tp P = I_{r_1}, \quad
        Q\tp Q = I_{r_2}, \quad S\tp S = I_{r_3}. 
      \end{equation}
using the Newton-Grassmann method (note that $P \in \RR^{k_1 \x
  r_1}$, $Q \in \RR^{k_2 \x r_2}$, and $S \in \RR^{k_3 \x r_3}$).    
Denote the core tensor after this computation by
$\widetilde{\cF}$.
Due to the fact that $\cF$  is a subtensor of
               ${\cC}$,   it follows  that 
      \begin{equation}
        \label{eq:Fnon-decr}
          \| \widetilde{\cF} \| \geq \| \cF \| ,       
       \end{equation}
    and evidently the Krylov-Schur  algorithm produces a
    non-decreasing sequence of 
    objective function values that is bounded above (by $\| \cA \|$).

    If $E_0$ is the local maximum point for
    \eqref{eq:maxB}, then  the G-gradient $\nabla_\cB(E_0)=0$,  and, by
    Proposition \ref{prop:gradient}, 
    $\nabla_\cC(\bar E_0)=0$,  where $\bar E_0$ corresponds to
    $E_0$. Therefore,  the 
    Newton-Grassmann method for \eqref{eq:maxrC}  will not
    advance, but give $ \widetilde{\cF} = \cF$. 

    On the other hand, if $E_0$ is not the local
    maximum, then $\nabla_\cB(E_0)$ and
    $\nabla_\cC(\bar E_0)$ are
    nonzero. Assume that we are close to a local maximum so that the
    G-Hessian for \eqref{eq:maxB} is negative definite. Then  the
    G-Hessian for \eqref{eq:maxC} is also negative definite, see
    Appendix\footnote{In the appendix we also take care of some
      Grassmann-technical details.} \ref{sec:app-hess},  and the projected
    maximization problem \eqref{eq:maxrC} is 
    locally convex. Therefore the Newton-Grass\-mann method will
    converge  to a solution that satisfies 
    $ \| \widetilde{\cF} \| > \| \cF \|$ \cite[Theorem 3.1.1]{flet87}.

    Thus, we have the following result.
    
    \begin{theorem}\label{theo:convergence}
    Assume that  the assumptions  \eqref{eq:stat-assump} hold, and that
    $(U_0,V_0,W_0)$ is close enough to a local maximum for
    \eqref{eq:maxrA}. Then  Algorithm \ref{alg:tks} will converge to that
    local maximum.       
    \end{theorem}

    The G-gradient is zero at the local maximum; thus the
    Krylov-Schur-like algorithm converges to a stationary point for the
    best rank-$(r_1,r_2,r_3)$ approximation problem.

      \section{Numerical tests}\label{sec:numerical test}

      In this section we will investigate the performance of the different
      block-Krylov-Schur methods by applying them to some test
      problems. As a comparison  we will use the HOOI method.  We here
      give a       brief description, for details, see
      \cite{de2000best}. 

 \subsection{Higher Order  Orthogonal Iteration}
 \label{sec:HOOI}

       Consider first      the non\-symmet\-ric case
      \eqref{eq:max}. HOOI is an  alternating  iterative method
      \cite{de2000best}, where in each 
      iteration three maximization problems are       solved.  In the
      first  maximization      we assume that $V$ and $W$ are given
      satisfying the      constraints, and maximize
      \begin{equation}\label{eq:HOOI}
        \max_{U^T U = I} \| \tmr[1]{\cC_1}{U}\|, \qquad \cC_1 =
        \tmr[2,3]{\cA}{V,W} \in \RR^{l \x r_2 \x r_3}.
      \end{equation}
      The solution of this problem is given by the first $r_1$ left
      singular vectors of the mode-1 unfolding $C^{(1)}$  of $\cC_1$, and that is
      taken as the new approximation $U$. Then in the second
      maximization $U$ and $W$ are considered as given and $V$ is
      determined, etc.

      The cost for computing
      the thin SVD  is $O(l (r_2 r_3)^2)$ (under the
      assumption \eqref{eq:rk-assumptions}). As this  computation is
      highly optimized, it is safe  to assume that for large and sparse
      tensors the computational cost in HOOI is dominated by
      the tensor-matrix multiplications $\tmr[2,3]{\cA}{V,W}$ (and
      corresponding in the other modes), and the reshaping of the
      tensor necessary for the multiplication. In addition, the
      computation of the G-gradient   involves  extra tensor-matrix 
      multiplications. Here it is more efficient to use global
      coordinates, cf. Section \ref{subsec:gradient}. In our experiments
      we computed the G-gradient       only every ten iterations. 

      For a (1,2)-symmetric tensor we use the HOOI method,
      where we have two maximizations in each step, one for $U$, with
      the previous value of $U$ in $\cC_1$, and the other
      for $W$, with $\cC_3=\tmr[1,2]{\cA}{U,U}$.  

      The HOOI iteration has linear
      convergence; its convergence properties are studied in
               \cite{xu15}. In general, alternating methods are  not
               guaranteed 
      to converge to a local  maximizer
      \cite{ruhe1980algorithms,powell1973search}. For some
      tensors, HOOI needs quite a few
      iterations before the convergence is stabilized to a constant linear
      rate. On the other hand,  HOOI has relatively fast convergence for
      several  well-conditioned problems. In our tests
      the HOOI method is initialized with random  matrices with
      orthonormal columns.

 For large tensors we use HOOI as a starting procedure for the
      inner Newton-Grassmann iterations to improve robustness. Note
      that here the tensor $\cC$  is much smaller of dimension
      $(k_1,k_1,k_3)$.

      \subsection{Experiments}
 \label{sec:exp}

 To investigate the performance of
 Krylov-Schur methods  we present in the following the
 results of  experiments on a few 
 $(1,2)$-symmetric tensors.
 In the  outer iteration, the  stopping criterion is 
 the relative norm of the G-gradient
 ${\norm{\nabla}}/{\norm\cF} \leq 
 10^{-13}$.

 In the inner iteration (Step (i) of algorithm \ref{alg:tks}) 
  the  best 
 rank-$(r_1,r_1,r_3)$ of the  tensor $\cC$ is computed 
 using  the Newton-Grassmann algorithm \cite{elden2009newton},
 initialized by  the 
 HOSVD of  $\cC$, truncated to rank $(r_1,r_1,r_3)$, followed by 5 HOOI
 iterations to improve robustness (recall that this problem is of
 dimension $(k_1,k_1,k_3)$).      The same
 stopping criterion as in the outer iteration was used. The typical
 iteration count for the inner Newton-Grassmann iteration was 4-5.
 
 In the  figures  the
 convergence history of the following four methods are  illustrated,
 where the  last three are named according to  the block Krylov-type method
               used in the outer  iterations: 
      \begin{enumerate}
      \item HOOI,
      \item  max-BKS($s$\,;\,$k_1$,$k_3$),
       \item BKS($s$,$p$\,;\,$k_1$,$k_3$),
         \item min-BKS($s$,$p$\,;\,$k_1$,$k_3$). 
         \end{enumerate}
      Here  $s$  denotes the  stage,  $k_1$ and $k_3$ the number of
      basis vectors in the $X$  and $Z$ basis, respectively, and $p$
      indicates that    the first $p$ vectors of each
      ON-block are used in BKS and min-BKS.

              The convergence results are presented in   figures, where
      the $y$-axis and $x$-axis represent   the relative norm of G-gradient
      and the number of outer iterations, respectively. For the larger
      examples we also plot the gradient against the execution time.
               
               In the examples we used rather small values for 
               $(r_1,r_1,r_3)$, like in some  real world applications
               \cite{liu2013multiview,persson2016maps}. 
               In two forthcoming papers \cite{eldehg20a,eldehg20c}
               on tensor  partitioning and                data science
               applications  
                we  compute rank-$(2,2,2)$ approximations.

                The experiments were performed using Matlab and the
                tensor toolbox \cite{bako:08} on a standard desktop
                computer with 8 GBytes RAM memory.  In all test cases the 
                memory was sufficient. 
                
      \subsubsection{Example 1. Synthetic  signal-plus-noise tensor}
 \label{sec:synth}
 For the first test we generated  synthetic $(1,2)$-symmetric tensors  with 
      specified low rank. Let $\cA=\cA_{signal}+\rho\cA_{noise}$,
               where $\cA_{signal}$ is a 
 signal tensor with low multilinear rank,
      $\cA_{noise}$ is a  noise tensor and $\rho$
      denotes the noise level. $\cA_{signal}$ was constructed as
       a   tensor of dimension $(r_1,r_1,r_3)$ with normally
      distributed $N(0,1)$ elements;  this was placed at  the top
       left corner of a zero tensor of  size $m \x m\x n$. The
      elements of the noise tensor were  chosen normally
               distributed $N(0,1)$,  and then the 
      tensor was made $(1,2)$-symmetric.  For testing purposes we 
      performed a random permutation such that the tensor remained
      (1,2)-symmetric. This tensor is dense. The purpose of this
      example is to demonstrate that the rate of convergence of all
      methods depends heavily on  the conditioning of the approximation
      problem (cf. \cite[Corollary 4.5]{elden2011perturbation} and the
      short statement in  Section \ref{sec:conv-KS}).
      
\begin{figure}[htb] 
   \centering
 \includegraphics[height=7cm,width=8cm]{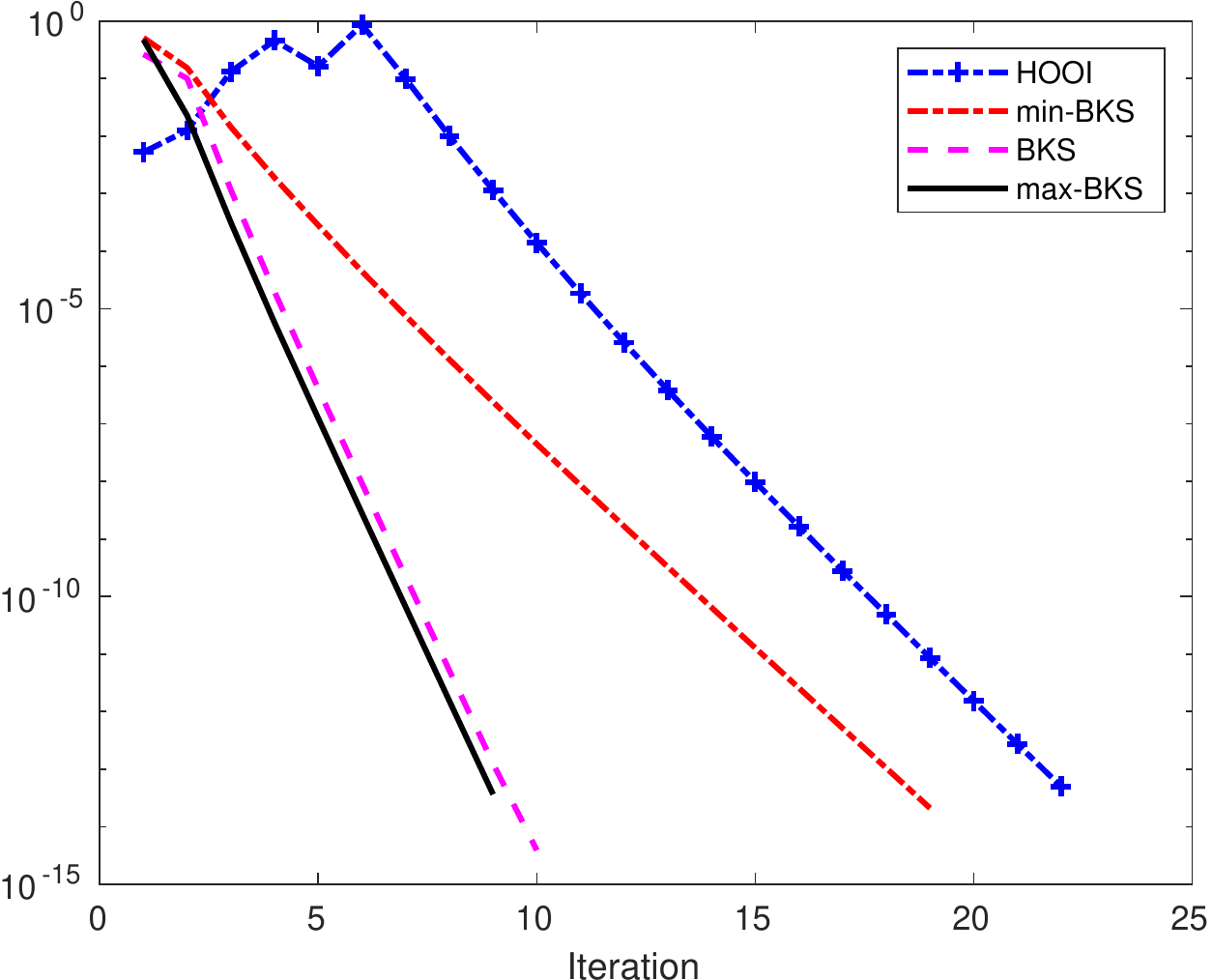}\\
 \includegraphics[height=7cm,width=8cm]{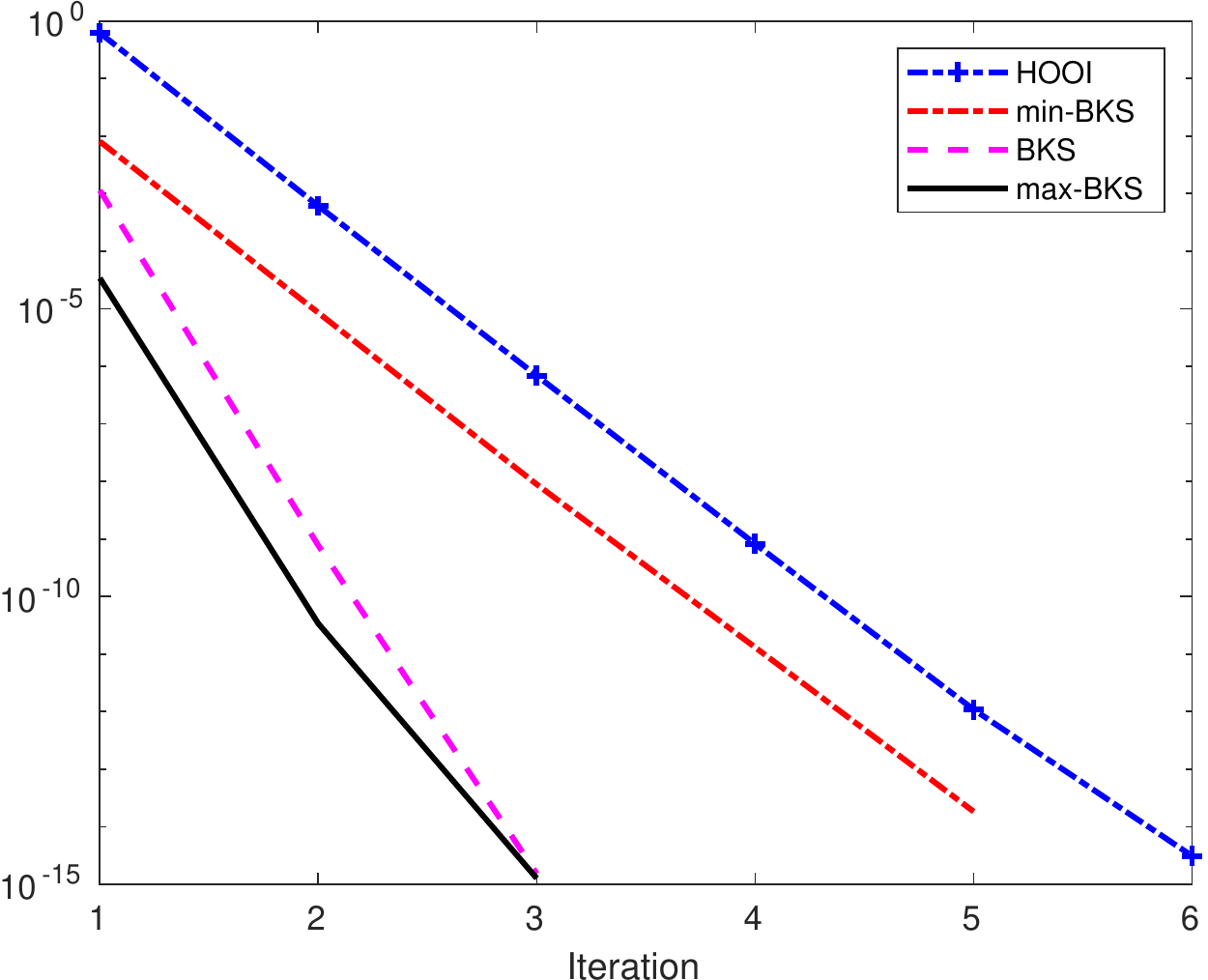}
   \caption{Convergence for  Example 1 when $(m,m,n)=(200,200,200)$
     and $(r_1,r_1,r_3)=(2,2,2)$. In the top 
      plot $\rho=10^{-2}$ and in the bottom  $\rho=10^{-4}$. In both
      cases the methods min-BKS(3,4;34,25), BKS(2,4;20,13), and
      max-BKS(2;32,23) were used. }
   \label{fig:synth}
 \end{figure}

      Figure \ref{fig:synth} illustrates the convergence results for a $200\x200
      \x 200$ tensor  approximated by one of rank-$(2,2,2)$, which is
      the correct rank of the signal tensor.
      The problem with $\rho=10^{-2}$ is more
      difficult than the other one, because the signal-to-noise ratio is
      smaller, making the approximation problem more ill-condi\-tioned.  This
      shows in the iteration count for all methods. See also Table
      \ref{tab:Ex1-S-values}, where the S-values are given. 
      The HOOI was
      considerably 
      more sensitive to the starting approximation than the other
      methods. For a few  starting values it converged  much more
      slowly than in the 
      figure.

 \begin{table}[htb]
   \centering
      \caption{Example 1. S-values for $\rho=10^{-2}$ (left) and $\rho=10^{-4}$
     (right). The column $k$ is the mode. }
   \begin{tabular}{|c|cc|c|}
     \hline
     $k$  & $s_1^{(k)}$ & $s_2^{(k)}$  & $s_3^{(k)}$  \\
     \hline
     1 & 1.95 & 1.33 & 0.11 \\
     3 & 2.07 & 1.13 & 0.14 \\
     \hline
   \end{tabular} \qquad
      \begin{tabular}{|c|cc|c|}
     \hline
     $k$  & $s_1^{(k)}$ & $s_2^{(k)}$  & $s_3^{(k)}$   \\
     \hline
     1 & 1.80 & 1.68 & 0.001 \\
     3 & 2.43 & 0.35 & 0.001 \\
     \hline
   \end{tabular}
   \label{tab:Ex1-S-values}
 \end{table}
 
 For this small problem it is possible to compute the solution
 directly using HOOI and the Newton-Grassmann method, without the
 Krylov-Schur approach. We compared that solution with the one
 obtained using the max-BK method and they agreed to
 within less than the magnitude of 
 the stopping criterion (the angles
   between the       subspaces spanned by the solution matrices were
   computed).



For the small problems in this example it is not meaningful to compare
the computational efficiency of the BKS methods to that of HOOI: due
to the simplicity of HOOI,  it comes out as a winner in the cases
when it converges.

       \subsubsection{Example 2. The Princeton  tensor}
 \label{sec:sparse-Prince}

The  Princeton tensor is created using the
      Facebook data from Princeton \cite{tkmp11}\footnote{The data
      can be downloaded from
      \url{https://archive.org/details/oxford-2005-facebook-matrix}.}.
      We constructed it   from a social network,  
      using a student/faculty flag  as third mode: the tensor element 
      $\cA(\lambda,\mu,\nu)=1$, if students $\lambda$ and $\mu$ are
      friends and one of them has flag $\nu$. After  zero slices in the third
      mode  have been removed, this is  a $6593\x6593\x29$
      tensor with $1.2 \cdot 10^6$ non zeros.

\begin{figure}[htb]  
   \begin{center} 
     \includegraphics[height=5cm,width=6cm]{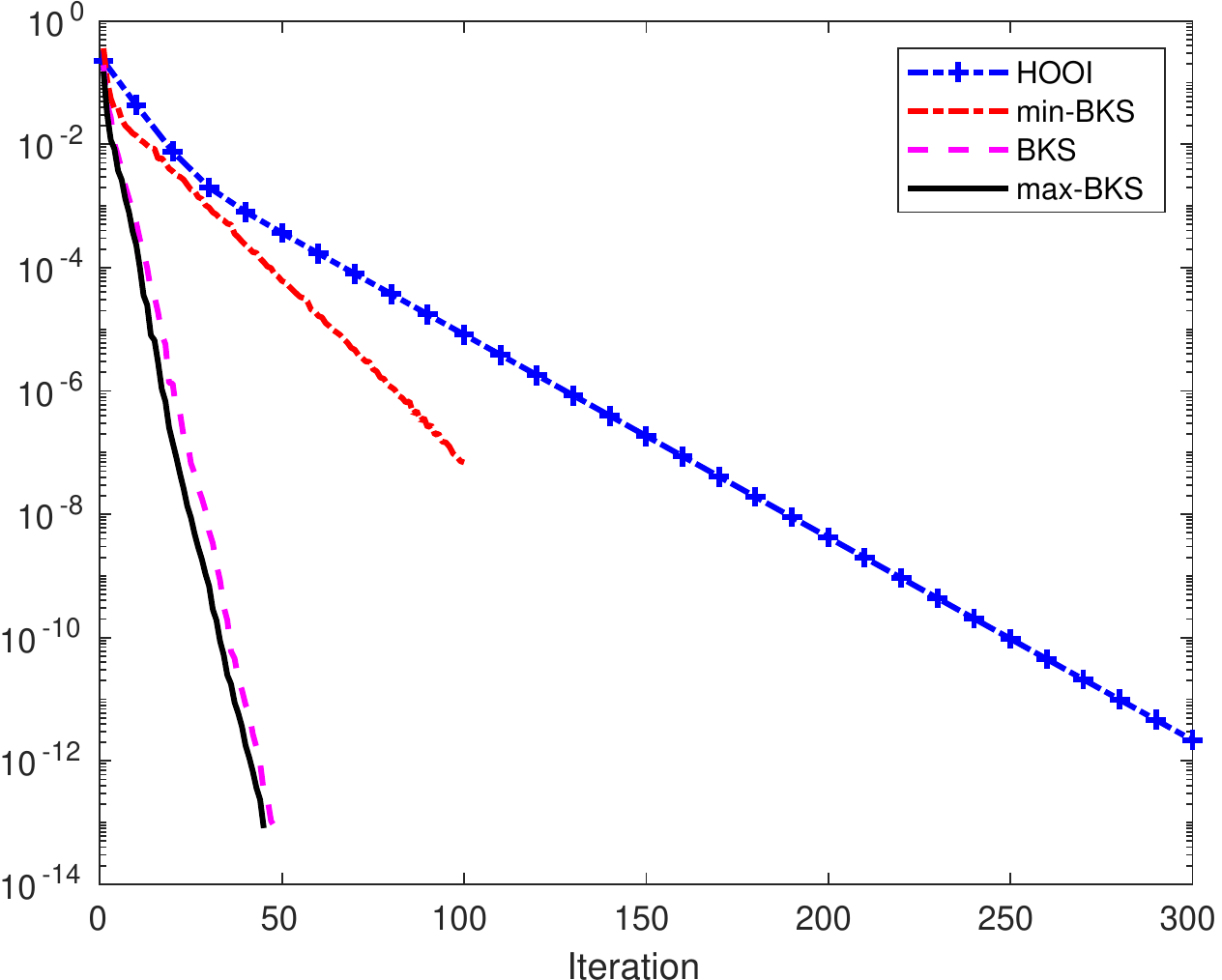}
      \includegraphics[height=5cm,width=6cm]{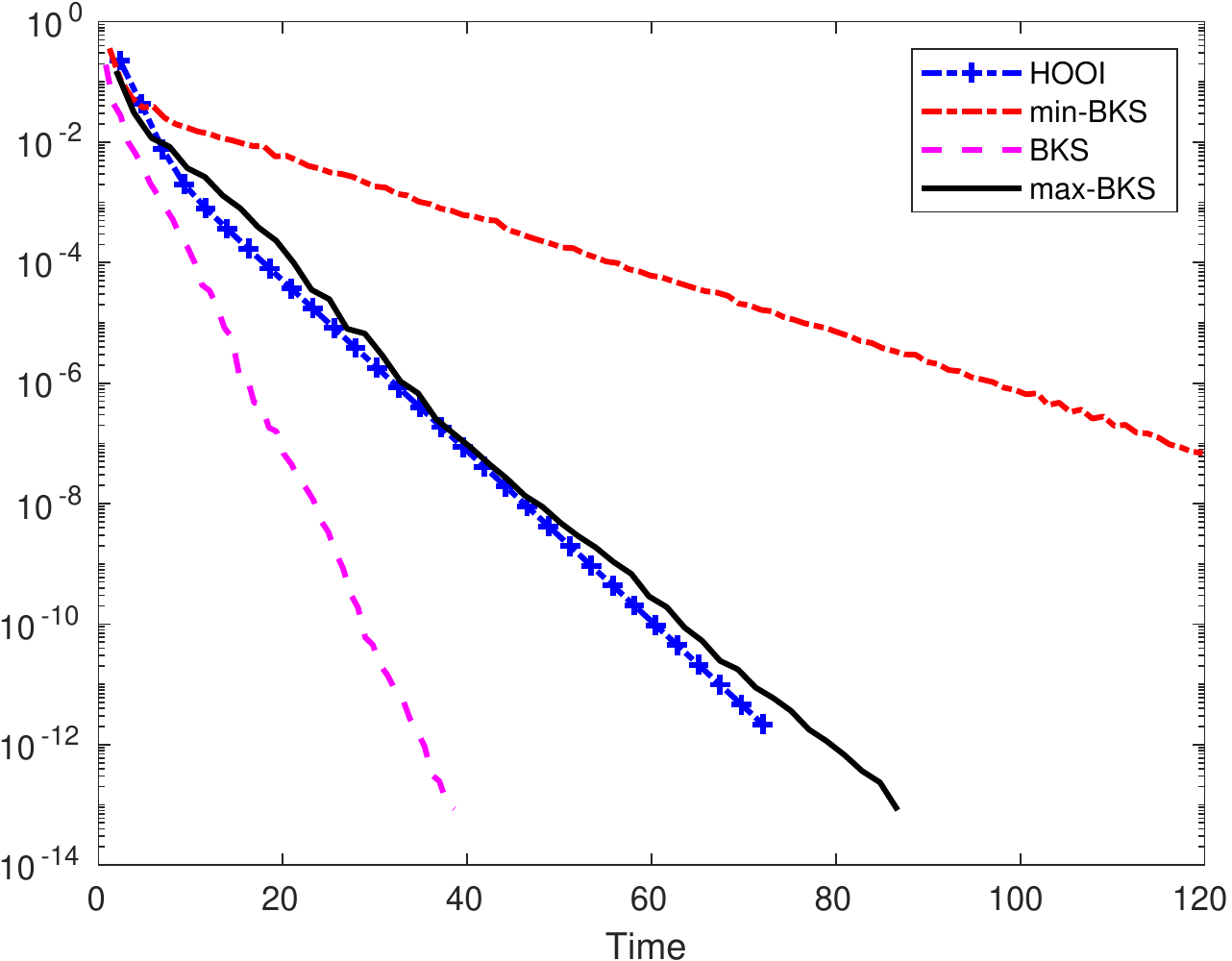}\\
  \caption{Convergence for  Example 2, the Princeton tensor with
    $(m,m,n)=(6593,6593,29)$ and     $(r_1,r_2,r_3)=(3,3,3)$.
    Convergence as function of 
    iterations (left),  as function of time (seconds), (right). 
 The methods min-BKS(3,5;62,29),
BKS(2,4;36,21), and max-BKS(1;12,9) were used.}
 \end{center}
   \label{fig:ex3princeton1} 
 \end{figure}

 Figure \ref{fig:ex3princeton1} shows the convergence  for the 
      Princeton tensor  approximated with a rank-$(2,2,2)$ and a
      $(4,4,4)$-rank tensor.  In both cases 
    the inner iterations   were
      initialized with truncated HOSVD followed by 5 HOOI iterations.
  
 The time measurements in Figure \ref{fig:ex3princeton1} are based on the
      Matlab functions \texttt{tic} and \texttt{toc}. A large
      proportion of the time in HOOI is the computation of the
      G-gradient (which was done every ten iterations).

 \begin{table}[htb]
   \centering
      \caption{Example 2, Princeton tensor. S-values for
        $(r_1,r_2,r_3)=(3,3,3)$.
        Note that the very small entry is approximately
     equal to zero in the floating point
     system.   \label{tab:Princeton-S-values}  }
   \begin{tabular}{|c|ccc|c|}
     \hline
     $k$  & $s_1^{(k)}$ & $s_2^{(k)}$  & $s_3^{(k)}$ & $s_4^{(k)}$  \\
     \hline
     1 & 300 & 193 & 187 & 47.3 \\
     3 & 390 & 185 & 80.2 & $2.4\cdot 10^{-14}$  \\
     \hline
   \end{tabular}
%

 %
 \end{table}

 In Table \ref{tab:Princeton-S-values} we give the S-values. The
 mode-3 entries are a strong indication that the mode-3 multilinear
 rank of the tensor is equal to 3. Any attempt to compute e.g. a
 rank-(4,4,4) approximation will suffer from the mode-3
 ill-conditioning and is likely to give incorrect results. However, a
 rank-(4,4,3) can be computed easily using BKS. Here the  convergence
 of HOOI was  very slow. 
 
       The number of iterations in the BKS method was
      rather insensitive to the choice of stage  and block
      size parameters $s$ and $p$. Thus it did not pay off to use a
      larger inner subproblem. Similarly, the use of max-BKS(2;111,29)
      gave relatively fast convergence in terms of the number of
      iterations, but  the extra information gained by using a
      large value of $k_1$ was not so substantial that it made up for
      the heavier computations.

      HOOI was sensitive to the choice of starting
      approximation. Often the convergence rate was considerably
      slower than in Figure \ref{fig:ex3princeton1}.

 \subsubsection{Example 3.  The Reuters tensor}
 \label{sec:sparse-Reuters}
    The Reuters tensor is a sparse tensor of dimensions  $13332\x
      13332 \x 66$ with $486,894$ nonzero elements. It is constructed
      from all stories released during 
      66 consecutive days by the news agency Reuters after the
      September 11, 2001, attack 
      \cite{batagelj2001density}. The vertices of the network are
      words. There is an edge between two words if they appear in the
      same text unit (sentence). The weight of an edge is its
      frequency.  

\begin{figure}[htb]  
   \begin{center}
     \includegraphics[height=5cm,width=6cm]{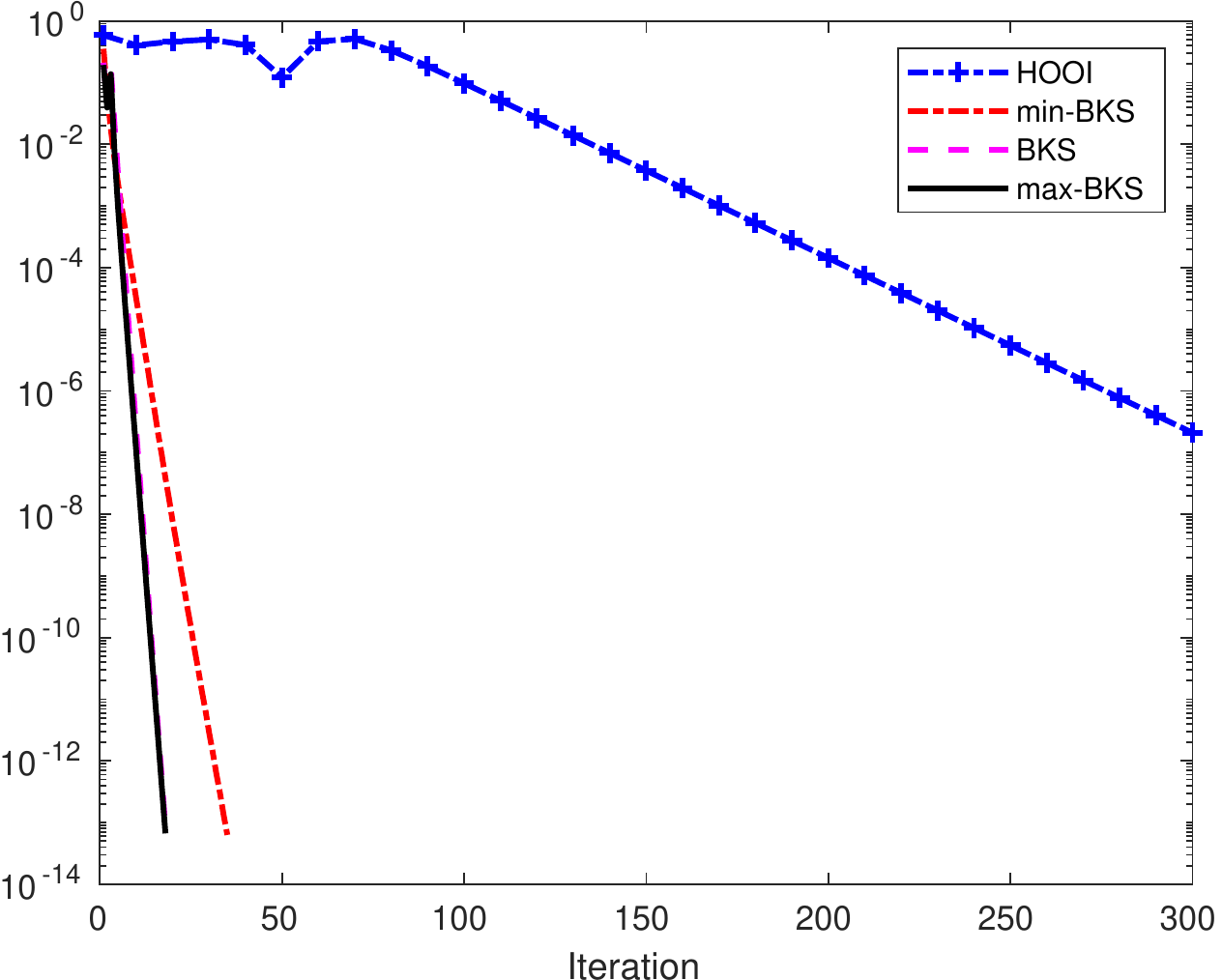}
      \includegraphics[height=5cm,width=6cm]{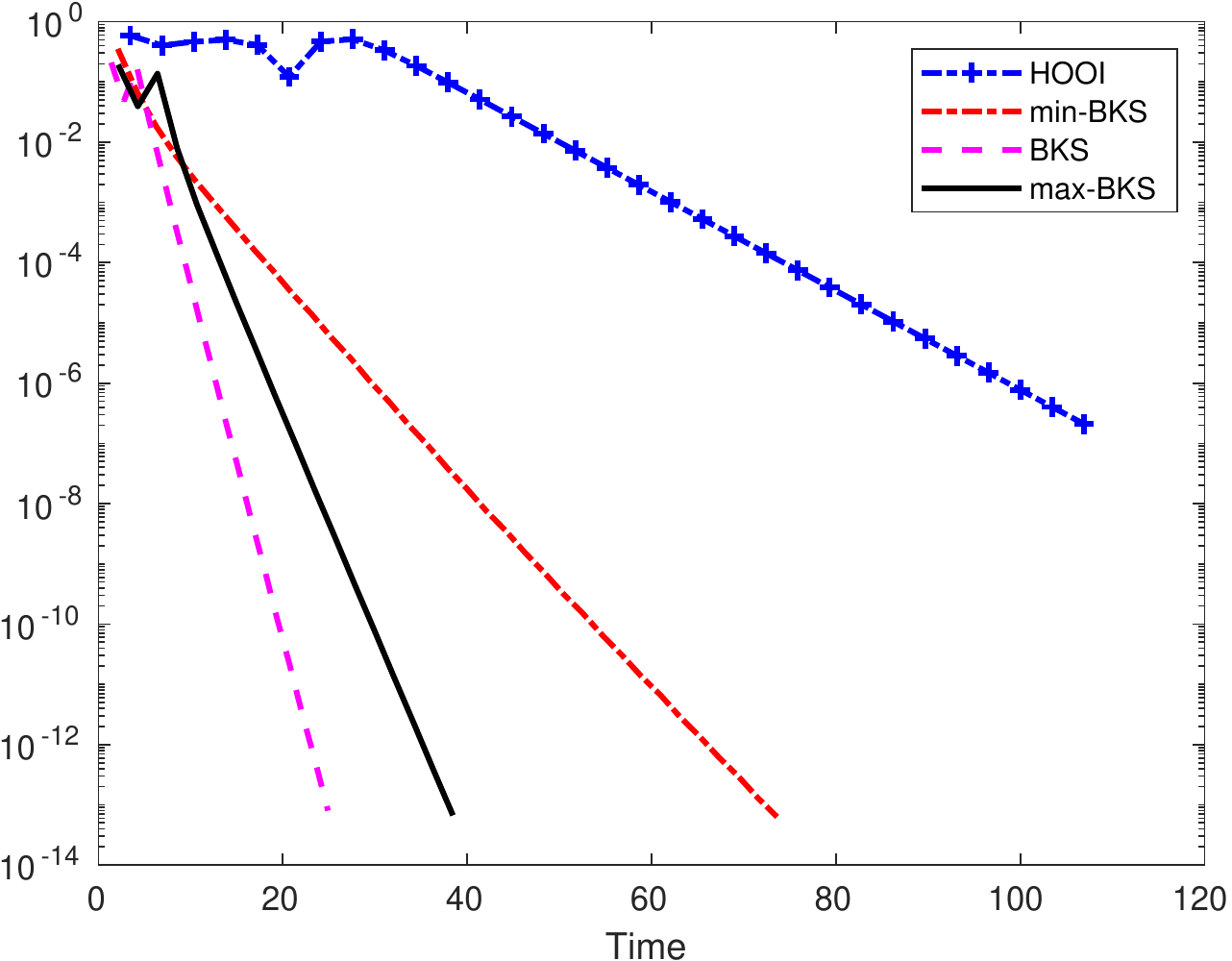}\\
      \includegraphics[height=5cm,width=6cm]{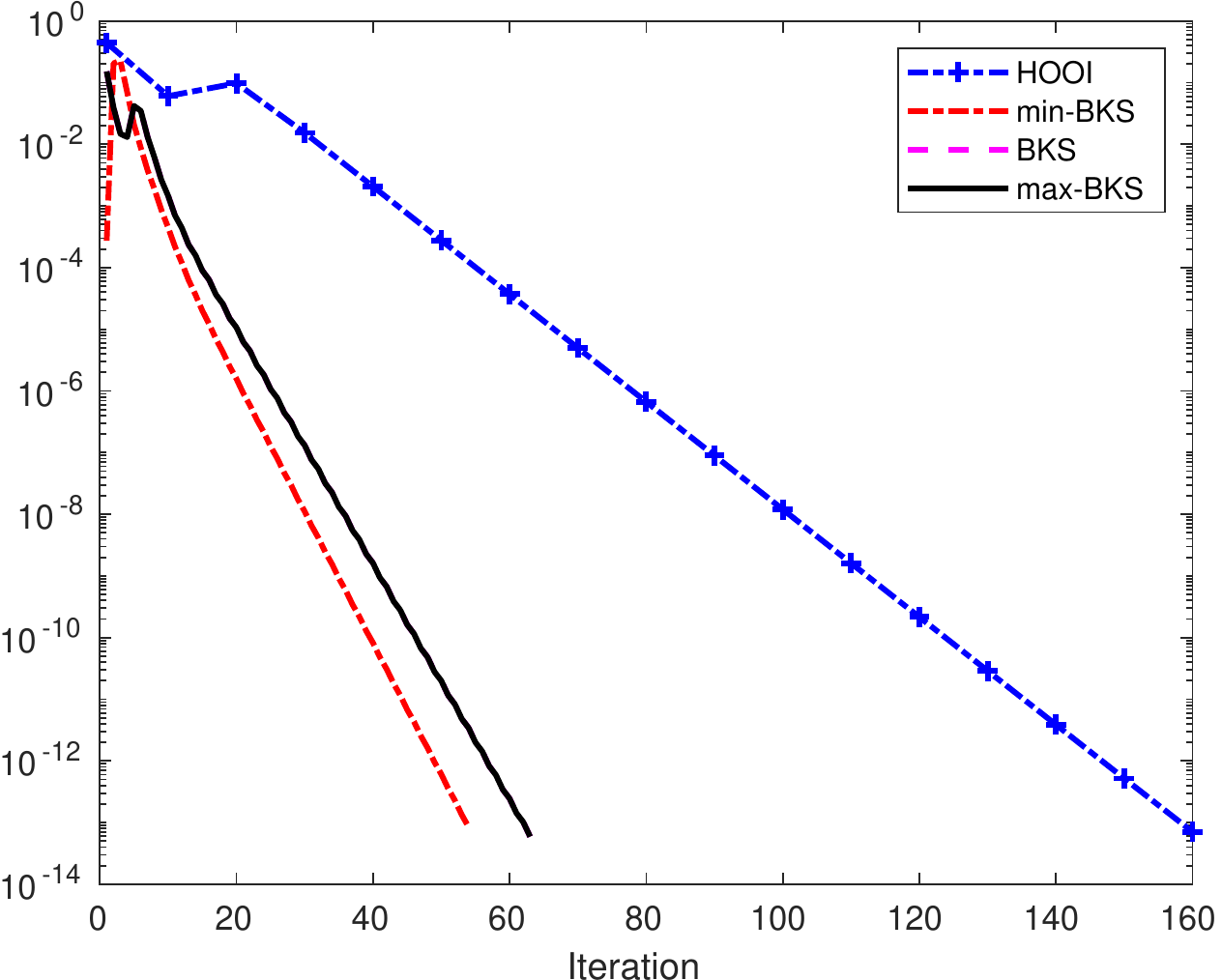}
       \includegraphics[height=5cm,width=6cm]{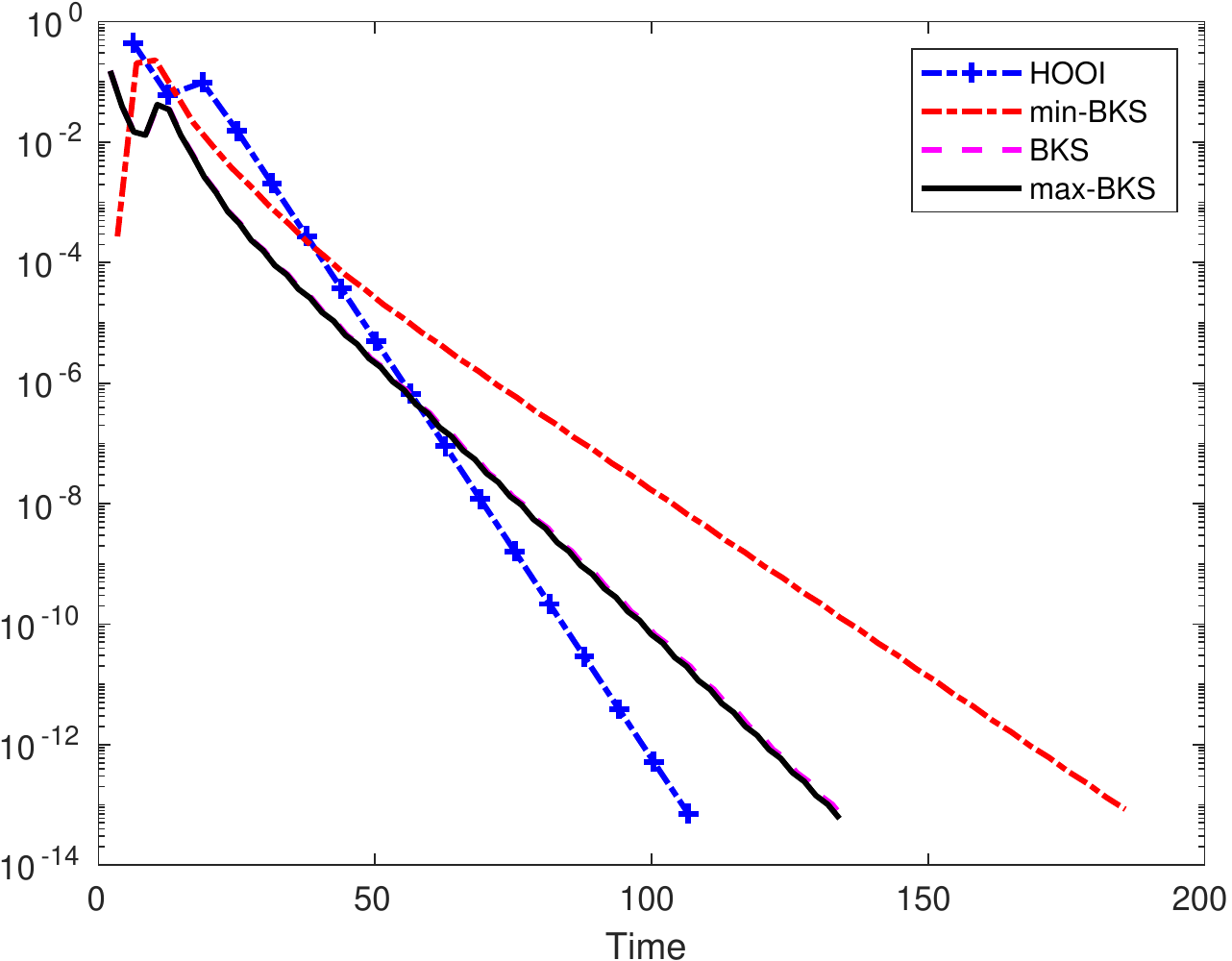}
  \caption{Convergence for Example 3, the Reuters tensor with
    $(m,m,n)=(13332,13332,66)$. 
     Top plot: $(r_1,r_2,r_3)=(2,2,2)$. The methods min-BKS(3,4;34,25),
     BKS(2,4;20,13)  and     max-BKS(2;32,23) were used. 
     Bottom plot: $(r_1,r_2,r_3)=(6,6,6)$. The methods min-BKS(2,4;58,37),
     BKS(1,6;42,27) and     max-BKS(1;42,27) were used (the last two
     are identical for these parameters). } 
 \end{center} 
   \label{fig:ex3reuters} 
 \end{figure}

 Figure \ref{fig:ex3reuters} shows the convergence results for the Reuters
      tensor, approximated with a rank-$(2,2,2)$ and  a rank-$(6,6,6)$
      tensor. In the second case, the inflexibility of the choice
      of $k_1$ and $k_3$ in max-BKS
      forced us to use stage 1, which led to slower
      convergence than with BKS and min-BKS. 
 
 The S-values are given in Table \ref{tab:Reuters-S-values}. It is
 seen that none of the  problems is particularly ill-conditioned.
 
 \begin{table}[htb]
   \centering
      \caption{Example 3 S-values for the Reuters tensor with
        $(r_1,r_2,r_3)=(2,2,2)$ (left) and  $(r_1,r_2,r_3)=(6,6,6)$ (right).}
   \begin{tabular}{|c|cc|c|}
     \hline
     $k$  & $s_1^{(k)}$ & $s_2^{(k)}$  & $s_3^{(k)}$  \\
     \hline
     1 & 350 & 207 & 75.2 \\
     3 & 397 & 88 & 13 \\
     \hline
   \end{tabular}
\qquad
     \begin{tabular}{|c|ccc|c|}
     \hline
     $k$  & $s_1^{(k)}$ & $\cdots$ & $s_6^{(k)}$  & $s_7^{(k)}$  \\
     \hline
     1 & 353 & $\cdots$ & 125 & 74.3 \\
     3 & 495 & $\cdots$ & 22.2 & 15.7 \\
     \hline
   \end{tabular}
   \label{tab:Reuters-S-values}
 \end{table}
 
\begin{figure}[htb]  
   \begin{center}
     \includegraphics[height=5cm,width=6cm]{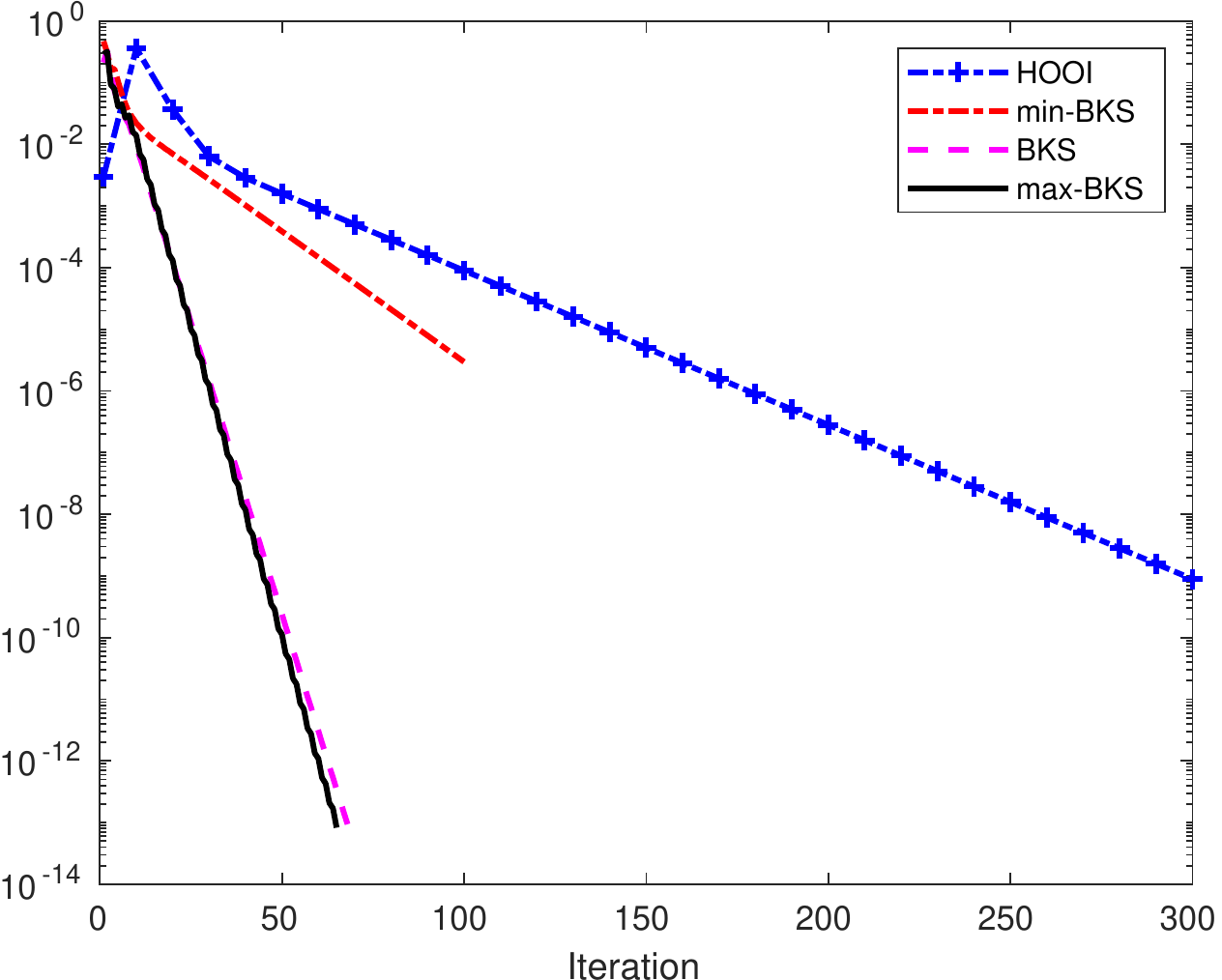}
      \includegraphics[height=5cm,width=6cm]{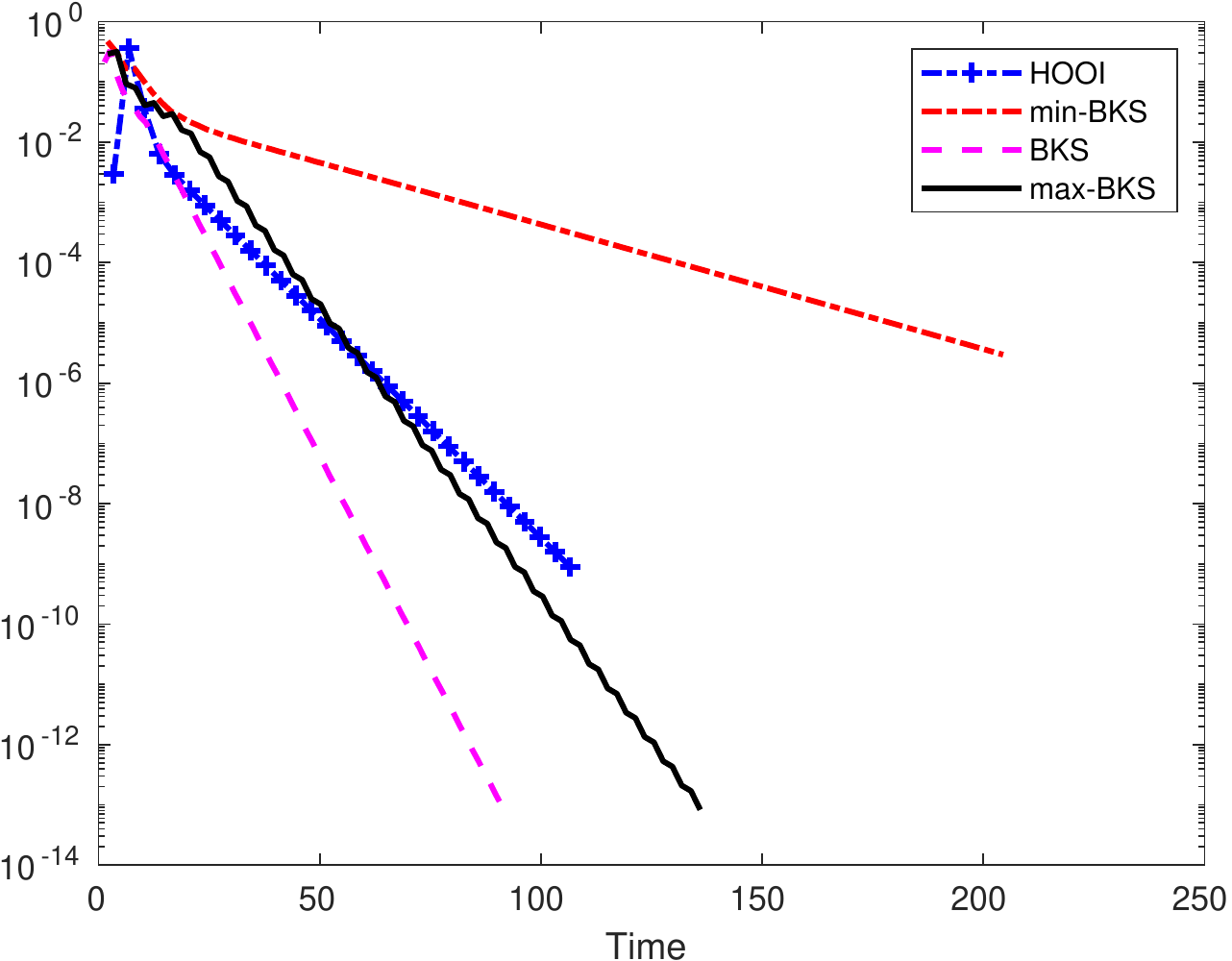}
  \caption{Convergence for Example 3, the scaled Reuters tensor, and
    $(r_1,r_2,r_3)=(2,2,2)$. We used min-BKS(3,4;34,25),
    BKS(2,4,32,23), and max-BKS(2;32,23).}
 \end{center} 
   \label{fig:ex3reutersscaled} 
 \end{figure}

 It is argued in \cite{eldehg20a} that in  cases when the 3-slices of
 a (1,2)-symmetric tensor are adjacency  matrices of graphs, then one
 should normalize the slices so that the largest eigenvalue of each
 slice becomes equal to 1. In that context  a rank-(2,2,2)
 approximation is computed. We ran a test with the normalized tensor
 and the same parameter values as in Figure \ref{fig:ex3reuters}. The
 results are shown in Figure \ref{fig:ex3reutersscaled}.
 
 The S-values are given in Table \ref{tab:ReutersScaledSValues}. They
 indicate that this problem is slightly more ill-conditioned than the
 unscaled one.

 \begin{table}[htb]
   \centering
      \caption{Example 3.  S-values for the scaled Reuters tensor with
        $(r_1,r_2,r_3)=(2,2,2)$.}
   \begin{tabular}{|c|cc|c|}
     \hline
     $k$  & $s_1^{(k)}$ & $s_2^{(k)}$  & $s_3^{(k)}$  \\
     \hline
     1 & 5.07 & 3.43  & 0.692 \\
     3 & 6.11  & 0.421 & 0.222 \\
     \hline
   \end{tabular}
   \label{tab:ReutersScaledSValues}
 \end{table}
 
\subsubsection{Example 4.  1998DARPA tensor}
\label{sec:1998}

\begin{figure}[htb]  
   \begin{center}
     \includegraphics[height=5cm,width=6cm]{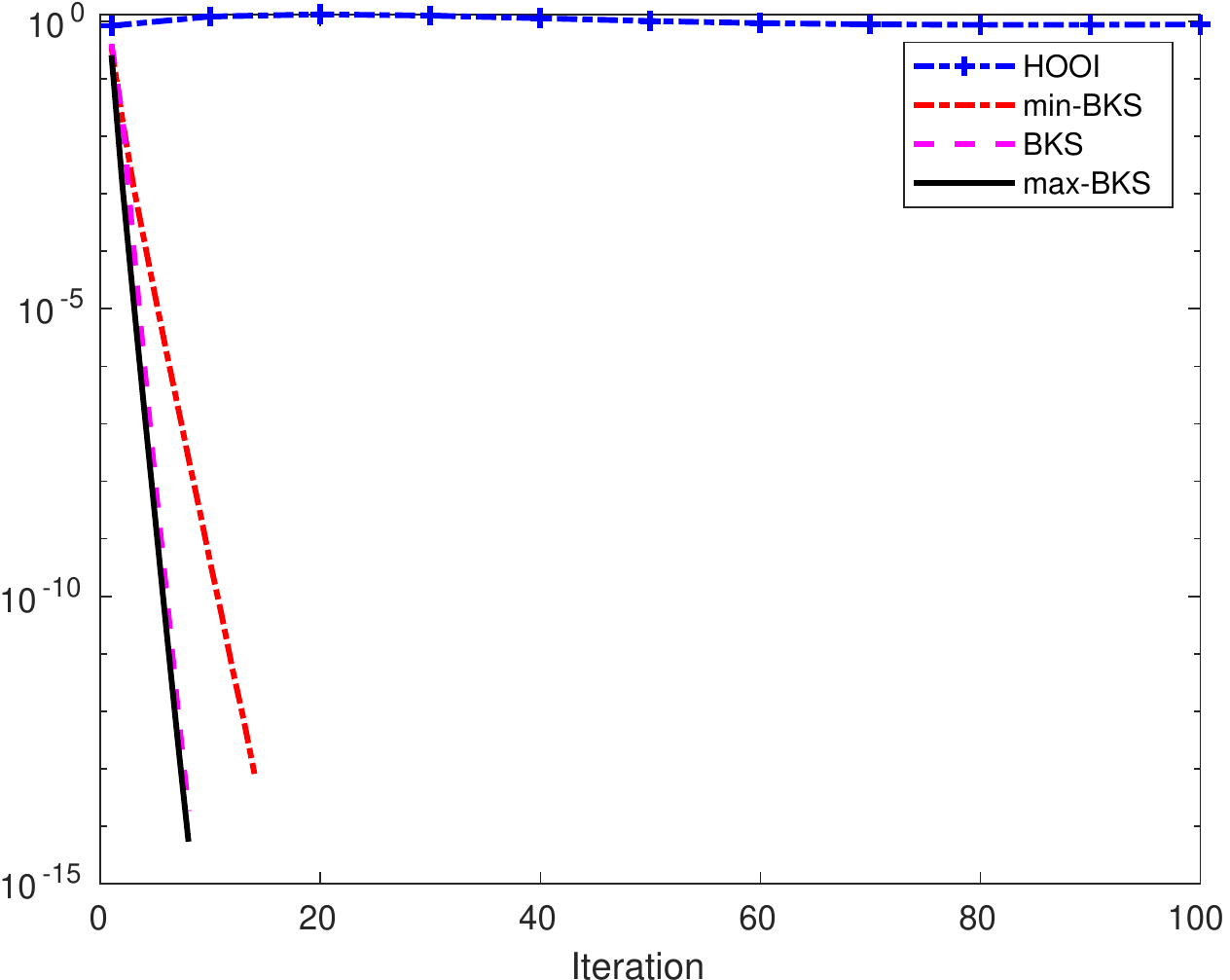}
      \includegraphics[height=5cm,width=6cm]{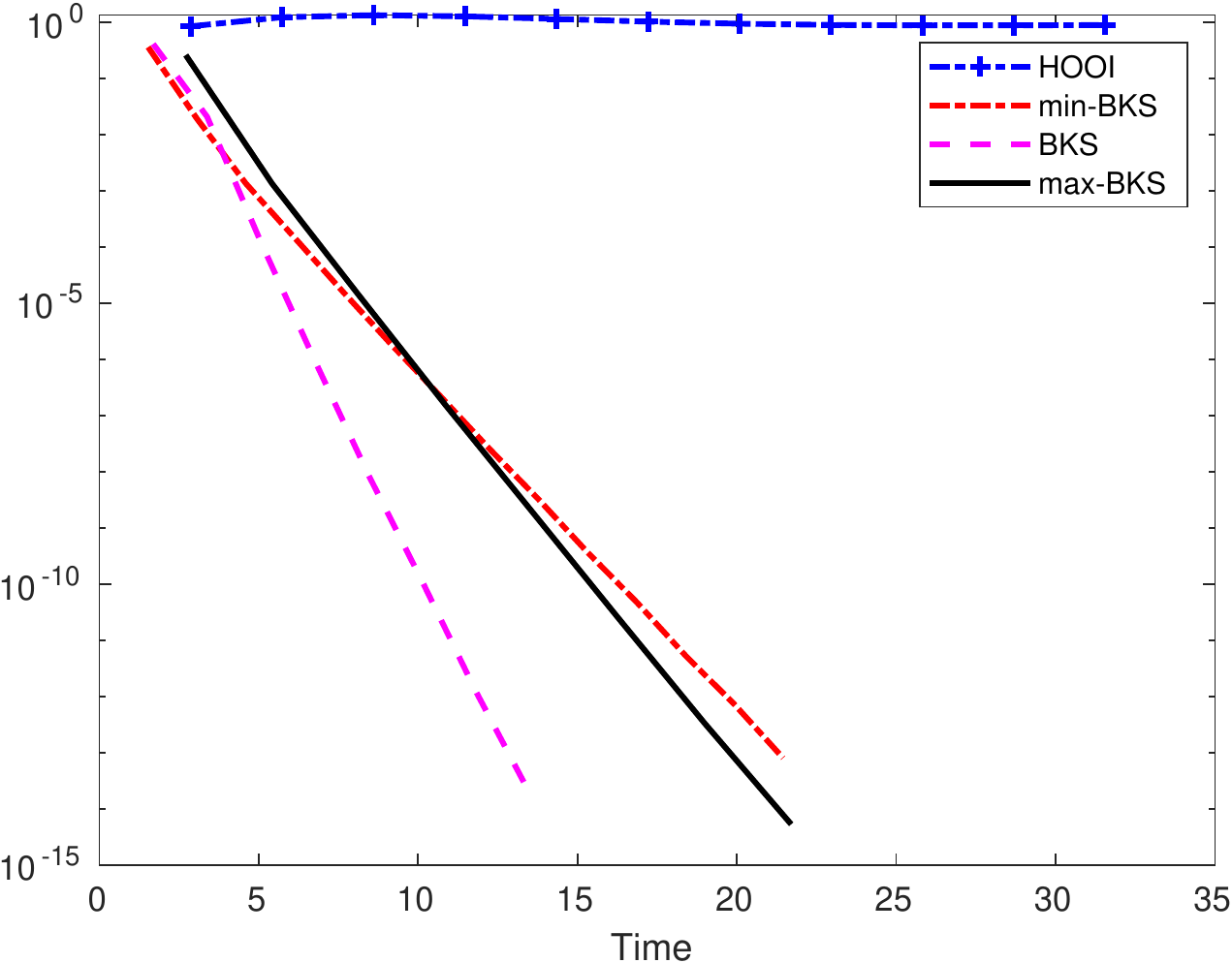}
  \caption{Convergence for Example 4, the 1998DARPA tensor with $(m,m,n)=(8991,8991,371)$, and
    $(r_1,r_2,r_3)=(2,2,2)$.  The methods
    min-BKS(2,4;18,15),  
     BKS(2,4;20,13)  and     max-BKS(2;32,23) were used.  } 
 \end{center} 
   \label{fig:ex41998DARPA} 
 \end{figure}

The following description is taken from \cite{eldehg20c}. In \cite{jpfsk16} network traffic logs are analyzed in order to
identify malicious attackers. The data are called the 1998 DARPA
Intrusion Detection Evaluation Dataset and were first published by the
Lincoln Laboratory at
MIT\footnote{\url{http://www.ll.mit.edu/r-d/datasets/1998-darpa-intrusion-detection-evaluation-dataset.}}. 
We downloaded the data set from
\url{https://datalab.snu.ac.kr/haten2/} in October 2018. The records
consist of (source IP, destination IP, port number, timestamp). In the
data file there are about 22000 different IP addresses. We chose the subset
of 8991 addresses that both sent and received messages. The time span
for the data is from June 1 1998 to July 18, and the number of
observations is about 23 million. We merged the data in time by 
collecting every 63999 consecutive observations into one bin. Finally
we symmetrized the tensor $\cA  \in \RR^{m \times m \times n}$,
where $m=8891$  and $n=371$,  
so that
\[
  a_{ijk}=
  \begin{cases}
    1 & \mbox{ if } i  \mbox{ communicated with } j \mbox{ in time
      slot } k \\
    0 & \mbox{ otherwise.}
  \end{cases}
  \]
In this example we did not normalize the slices of the
tensor: The 3-slices are extremely  sparse, and normalization makes
  the rank-(2,2,2) problem so ill-conditioned that none of the
  algorithms  converged. Instead we scaled the slices to have
  Frobenius norm   equal to 1. The convergence history is shown in
  Figure \ref{fig:ex41998DARPA}. 
 
 The HOOI method was sensitive to  the (random) starting
 approximations. It did happen that the method converged  rapidly,
 but in many cases convergence was extremely slow.

 The S-values are given in Table \ref{tab:ex41998DARPA}. The problem
 is well-conditioned.

 \begin{table}[htb]
   \centering
      \caption{Example 4.  S-values for the 1998DARPA tensor with
        $(r_1,r_2,r_3)=(2,2,2)$.}
   \begin{tabular}{|c|cc|c|}
     \hline
     $k$  & $s_1^{(k)}$ & $s_2^{(k)}$  & $s_3^{(k)}$  \\
     \hline
     1 & 6.35 & 6.12  & 0.948 \\
     3 & 8.82  & 0.257 & 0.0995 \\
     \hline
   \end{tabular}
   \label{tab:ex41998DARPA}
 \end{table}
 
\subsubsection{Discussion of experiments}
\label{sec:discussion}

Profiling tests of the BKS and HOOI methods for the Reuters example show
that most of the computational work is  done in the reshaping of
tensors, and in  tensor-matrix multiplications. For small values of
the rank $(r_1,r_1,r_3)$, the number of stages and block size in
the BKS methods,   the time work for the       dense tensor and
matrix operations in the  
inner iterations in BKS and the SVD's in HOOI is relatively small. A
considerable proportion of the work in HOOI is the computation
of the G-gradient; we reduced that by computing it only every ten
iterations.  The data shuffling and reshaping must be 
done irrespectively of which programming language is used. It is
reasonable to assume that the implementation made in the sparse
tensor toolbox \cite{bako:08} is efficient. 
Therefore it is makes sense to measure efficiency by comparing
the Matlab execution times (by \texttt{tic} and \texttt{toc)}
of the methods.
      
Our tests with (1,2)-symmetric tensors indicate that all
methods considered here       converge fast for very well-conditioned
problems. However, the convergence behavior of HOOI was less
predictable: sometimes it converged very slowly also for
well-conditioned, large problems.  Consistently the min-BKS method
converged much more slowly  
than the other two Krylov-Schur like methods. The max-BKS method 
suffered from its       inflexibility in the choice of $k_1$ and
$k_3$, especially with $r_1$ and $r_3$ somewhat larger. 

The design of BKS is to some extent based on heuristics and
      numerical experiments. A comparison of  BKS and max-BKS shows
      that the choice of blocks of $p$ vectors, for $p$  rather small,
      in the Krylov steps, does not
      substantially impede the convergence rate. 
      Using fewer blocks, in the sense of using only $p$ vectors from
      the ``diagonal'' blocks in the diagram in Table \ref{table:min-BK}, as 
      in min-BKS, leads to slower convergence. Thus BKS seems to be a
      reasonable compromise.  Based on the experience
      presented in  this paper and \cite{eldehg20a,eldehg20c} it seems
       clear that for large and sparse tensors the BKS method is in
       general  more robust and efficient       than HOOI.  

In the BKS method the parameters $s$ and  $p$ (which give $k_1$
and $k_3$) could be chosen rather small, typically 2 and 4,
respectively. Using larger values did not pay off.

 \section{Conclusions and future work}\label{sec:conclusions}
 
We have  generalized block Krylov-Schur methods for matrices to tensors 
      and demonstrated that the new method  can be used for computing
      best rank-$(r_1,r_2,r_3)$  
      approximations of large and sparse tensors. The  BKS method  is
      shown to be flexible and has  best convergence properties. 

      The purpose of this paper has been to show that the
      block-Krylov-Schur method  is a viable approach. 
      It may be possible to analyze BKS methods in depth,
      theoretically and by experiments,  and optimize
      the method further, for instance with regard to the choice of blocks in
      the tables defining the method.
      
      Since we are interested in very low rank approximation  of
      (1,2)-symmetric tensors for applications such as those in
      \cite{eldehg20a,eldehg20c},  the development of block
      Krylov-Schur type
      methods was done mainly with such applications in mind. More
      work is needed to 
      investigate the application of BKS methods for nonsymmetric
      tensors.
 
      The detailed implementation of block Krylov-Schur methods for
      matrices is rather technical, see e.g. \cite{zhsa08}. The
      generalization to tensors might improve the convergence
      properties for ill-conditioned problems. However, this is beyond
      the scope of the present paper, and may be a topic for future
      research. 

      It appears to be straightforward to generalize the results to
      tensors of order larger than 3. We are planning to do research
      in this direction in the future. 
      
\section{Acknowledgments} 
This work was done when the second author visited the Department of
      Mathematics, Link\"oping University.

 \bibliographystyle{siamplain}
       \bibliography{/media/lars/ExtHard/WORK/forskning/BIBLIOGRAPHIES/general,references,/media/lars/ExtHard/WORK/forskning/BIBLIOGRAPHIES/LE-papers}

      \bigskip

      \appendix
      
\section{The Grassmann Hessian}
\label{sec:app-hess}
Let $X \in \RR^{l \x r}$, be a matrix with orthonormal columns,  $X\tp
X=I_r$. We will let it represent an  entire subspace, i.e. the
equivalence class, 
\[
  [X] = \{ X Q \;|\; Q \in \RR^{r \x r}, Q\tp Q=I_r\}.
\]
 For convenience we will say that $X \in \Gr(l,r)$, the Grassmann
 manifold (of equivalence classes).

 Define the product manifold
 \[
   \Gr^3 = ( \Gr(l,r_1), \Gr(m,r_2), \Gr(n,r_3)),
 \]
 and, for given integers satisfying $r_1 < k_1 < l$, $r_2 < k_2 < m$,
 and $r_3 < k_3 < n$,
 \[
   \Gr^3_k = ( \Gr(k_1,r_1), \Gr(k_2,r_2), \Gr(k_3,r_3)).
 \]
 The following is a submanifold of $\Gr^3$:
 \[
   \Gr_s^3 = \{ (X,Y,Z) =
     \left(
       \begin{pmatrix}
         U \\ 0
       \end{pmatrix},
       \begin{pmatrix}
         V \\0
       \end{pmatrix},
       \begin{pmatrix}
         W \\ 0
       \end{pmatrix}
     \right) \in \Gr^3
     \;|\; (U,V,W) \in \Gr_k^3
   \}.
 \]
 Let $(X_0,Y_0,Z_0) \in \Gr^3$, and let three matrices  $X_1 \in \RR^{l \x
   (k_1-r_1)}$,  $Y_1 \in \RR^{m \x (k_2-r_2)} $, and $Z_1 \in \RR^{n
   \x (k_3-r_3)}$ be given, such that 
\[
\bar X =   [X_0 \, X_1], \qquad
\bar Y = [Y_0 \, Y_1], \qquad
\bar Z = [Z_0 \, Z_1],
\]
all have orthonormal columns.

For a given matrix $P \in \Gr(l,r)$ we let $P_\perp \in \Gr(l,l-r)$
 be such that $(P \, P_\perp)$ is an orthogonal matrix. It can be
 shown \cite[Section 2.5]{eas98} that $P_\perp$ is a matrix of basis vectors in the
 tangent space of $\Gr(l,r)$ at the point $P$. 

 In the maximization problem for $\| \tmr{\cA}{X,Y,Z}\|^2$ on the
 Grassmann product manifold $\Gr^3$, we now make a change of variables
 by defining
 \[
   \cB = \tmr{\cA}{(\bar X \, \bar X_\perp), (\bar Y \, \bar Y_\perp),
     (\bar Z \, \bar Z_\perp)},
 \]
 and further 
 \[
   \cC = \tmr{\cA}{\bar X , \bar Y ,
     \bar Z }.
 \]
 Clearly $\cC$ is a leading subtensor of $\cB$, see Figure \ref{fig:CB}.  After this change of
 variables the point $(X_0,Y_0,Z_0)$ is represented by
 \[
   \Gr^3 \ni E_0 = \left(
     \begin{pmatrix}
       I_{r_1} \\ 0
     \end{pmatrix},
     \begin{pmatrix}
       I_{r_2} \\ 0
     \end{pmatrix},
     \begin{pmatrix}
       I_{r_3} \\ 0
     \end{pmatrix}
   \right), 
 \]
 and the bases for the tangent space of $\Gr^3$  at $E_0$ are 
 \[
    (E_0)_\perp = \left(
     \begin{pmatrix}
       0 \\ I_{l-r_1} 
     \end{pmatrix},
     \begin{pmatrix}
       0 \\ I_{m-r_2} 
     \end{pmatrix},
     \begin{pmatrix}
       0 \\ I_{n-r_3} 
     \end{pmatrix}
   \right). 
 \]
 The bases  for the tangent space of the submanifold $\Gr_s^3$ at $E_0$ 
 are given by
  \[
   \left(
     \begin{pmatrix}
       0 \\ I_{k_1-r_1} \\ 0
     \end{pmatrix},
     \begin{pmatrix}
       0 \\ I_{k_2-r_2} \\ 0
     \end{pmatrix},
     \begin{pmatrix}
       0 \\ I_{k_3-r_3} \\0 
     \end{pmatrix}
   \right), 
 \]
 where the top zeros are in $\RR^{r_i \x r_i}$, $i=1,2,3$,  and the bottom in
 $\RR^{(l-k_1) \x (k_1-r_1)}$,  $\RR^{(m-k_2) \x (k_2-r_2)}$, and
      $\RR^{(n-k_3) \x (k_3-r_3)}$, respectively. Clearly, the tangent
      space of $\Gr^3_s$ is a subspace of the tangent space of
      $\Gr^3$.
      
 Define the functions
 \begin{align*}
   f(X,Y,Z) &= \| \tmr{\cB}{X,Y,Z} \|^2, \qquad (X,Y,Z) \in \Gr^3,\\
   g(U,V,W) &= \| \tmr{\cC}{U,V,W} \|^2, \qquad (U,V,W) \in \Gr_k^3
 \end{align*}
 The subtensor property implies that for $(U,V,W) \in \Gr_k^3$,
 \begin{equation}
   \label{eq:f=g}
   g(U,V,W) = f(X,Y,Z), \qquad
   (X,Y,Z) =
     \left(
       \begin{pmatrix}
         U \\ 0
       \end{pmatrix},
       \begin{pmatrix}
         V \\0
       \end{pmatrix},
       \begin{pmatrix}
         W \\ 0
       \end{pmatrix}
     \right). 
 \end{equation}
 \begin{proposition}
   \label{prop:H-pos-def}
   Assume that the Grassmann Hessian of $f$ is positive definite on
   the tangent space of $\Gr^3$ at $E_0$. Then the Grassmann Hessian
   of $g$ is positive definite on the tangent space of $\Gr_k^3$ at
   the point
   \[     
      E_{0k}  = \left(
     \begin{pmatrix}
       I_{r_1} \\ 0
     \end{pmatrix},
     \begin{pmatrix}
       I_{r_2} \\ 0
     \end{pmatrix},
     \begin{pmatrix}
       I_{r_3} \\ 0
     \end{pmatrix}
   \right) \in \Gr^3_k.
 \]
\end{proposition}

\begin{proof}
  As the tangent space in $\Gr^3_s$ at $E_0$ is a subspace of the
  tangent space in $\Gr^3$, the Hessian of $f$  must be positive
  definite at $E_0$ in $\Gr_s^3$. Therefore, due to \eqref{eq:f=g},
  the geometric properties of $g$  are the same as those of $f$, and
  the Hessian of $g$ is positive definite in $\Gr_k^3$ at $E_{0k}$.  
\end{proof}

\medskip

\section{Implementation of the block-Krylov step}
\label{sec:block-krylov}
The steps (ii)-(iv) in Algorithm \ref{alg:block-krylov-step} are
      written in tensor form to emphasize the equivalence to
      Gram-Schmidt orthogonalization. As we remarked in Section
      \ref{sec:block-krylov-methods}, the reorganization of data from
      tensor to matrix form before performing tensor-matrix
      multiplication  is costly. Therefore, we keep the result of
      step (i) in matrix form and directly orthogonalize it to the
      previous vectors by performing a QR decomposition. Thereby we
      also avoid performing reorthogonalization, which might be
      necessary if we use the Gram-Schmidt method.

 \end{document}